\documentclass[MSNbibl,number,citesort,dvips]{arxbj}
\usepackage{dcolumn}
\usepackage{graphicx}

\aid{0}
\volume{20}
\issue{1}
\pubyear{2014}
\firstpage{78}
\lastpage{108}
\doi{10.3150/12-BEJ477} %kopijuoti is PTS

\makeatletter
\newcolumntype{d}[1]{D{.}{.}{#1}}

\newcommand{\eqref}[1]{(\ref{#1})}
\newcommand{\field}[1]{\mathbb{#1}}
\newcommand{\R}{\field{R}}
\newcommand{\p}{\field{P}}
\newcommand{\N}{\field{N}}
\newcommand{\Z}{\field{Z}}
\newcommand{\E}{\field{E}}
\newcommand{\FF}{\mathcal{F}}
\newcommand{\GG}{\mathcal{G}}
\newcommand{\PP}{\mathcal{P}}
\newcommand{\RR}{\mathcal{R}}
\newcommand{\betav}{\bolds{\beta}}
\newcommand{\etav}{\bolds{\eta}}
\newcommand{\tr}{\operatorname{tr}}

\def\argmin{\mathop{\operatorname{argmin}}}
\newtheorem{theorem}{Theorem}
\newremark{example}{Example}
\newremark{remark}{Remark}
\newtheorem{lemma}{Lemma}
\newproclaim{definition}{Definition}

\newproclaim{condition}{Condition}
\newtheorem{proposition}{Proposition}

\def\RSS{\mathrm{RSS}}
\def\var{\operatorname{Var}}
\def\barid{\overline{ID}}
\makeatother

\begin{document}
\begin{frontmatter}

\title{Nonparametric specification for non-stationary time series regression}
\runtitle{Specification for non-stationary time series}

\begin{aug}
%%%% inicialai - be tarpu
\author{\fnms{Zhou} \snm{Zhou}\corref{}\ead[label=e1]{zhou@utstat.toronto.edu}}% \and
\runauthor{Z. Zhou} %% auto
\address{Department of Statistics, University of Toronto, 100 St. George Street, Toronto, Ontario,
M5S 3G3 Canada. \printead{e1}}
\end{aug}

% HISTORY:
\received{\smonth{1} \syear{2012}}
\revised{\smonth{8} \syear{2012}}

% ABSTRACT
%
\begin{abstract}
We investigate the behavior of the Generalized
Likelihood Ratio Test (GLRT) (Fan, Zhang and Zhang [\textit{Ann. Statist.}
\textbf{29} (2001) 153--193]) for time
varying coefficient models where the regressors and errors are
non-stationary time
series and can be cross correlated. It is found that the
GLRT retains the minimax rate of local alternative detection under
weak dependence and non-stationarity. However, in general, the Wilks
phenomenon as well as the classic residual bootstrap are sensitive to either
conditional heteroscedasticity of the errors,
non-stationarity or temporal dependence. An averaged test is suggested to
alleviate the sensitivity of the test to the choice of bandwidth and is
shown to be more powerful
than tests based on a single bandwidth. An alternative wild bootstrap
method is proposed and
shown to be consistent when making inference of time varying
coefficient models for non-stationary time series.
\end{abstract}

% KEYWORDS
% visi is mazosios raides ir pagal abecele
%

\end{frontmatter}

%s1 #&#
\section{Introduction}
\label{secintr}

Specification tests are important in many nonparametric settings.
Generally, one is interested in testing whether certain nonparametric
components are significant, or whether they have a more parsimonious
and efficient parametric representation. In the time series context,
there is a large literature devoting to the latter topic, see for
instance Hjellvik \textit{et al.} \cite{HjeYaoTjs98}, Fan and Li \cite
{FanLi99}, Dette and
Spreckelsen \cite{DetSpr03,DetSpr04}, An and Cheng \cite{AnChe91} and
Paparoditis \cite{Pap00},
among others. Many of the previous results perform specification for
stationary time series.

%Bickel and Rosenblatt (1973), Murphy (1993), Eubank and Speckman
%(1993), H\"{a}rdle and Mammen (1993), Ingster (1993), Inglot and
%Ledwina (1996), Hart (1997), Stute (1997), Xia (1998), Kauermann and
%Tutz (1999) and Horowitz and Spokoiny (2001), among others. Many of
%the previous results perform specification for independent samples.

The purpose of the paper is to develop specification tests for
nonparametric regression of non-stationary time series. Specifically,
consider the following time-varying coefficient model:
%
%e1 #&#
%
\begin{equation}
\label{eqmodel} y_i = \mathbf{x}_i^\top
\betav(t_i) + \varepsilon_i,\qquad  i=1,\ldots, n,
\end{equation}
where $t_i=i/n$, $\mathbf{x}_i = (x_{i1}, x_{i2}, \ldots,
x_{ip})^{\top}$
are $p
\times1$ dimensional time series of regressors or predictors,
$\varepsilon_i$ are
error series satisfying $\E(\varepsilon_i | \mathbf{x}_i) = 0$. Here
$^{\top}$ denotes matrix or vector
transpose. The processes $\{\mathbf{x}_i\}$ and $\{\varepsilon_i\}$ are
allowed to be non-stationary and can be cross correlated. We assume
that the regression parameters $\betav(\cdot) := (\beta_1(\cdot),
\ldots
, \beta_p(\cdot)
)^{\top}$ is a smooth function on $[0, 1]$. Nonparametric specification
of model \eqref{eqmodel} boils down to testing whether $\betav(\cdot)$
or a component of it has a certain parametric representation.

Due to their flexility and interpretability in investigating shifting
association between the response and predictors over time, model \eqref
{eqmodel} and its stochastic coefficient version have attracted
considerable attention in various fields. See, for instance, Orbe
\textit{et al.}
\cite{OrbFerRod05,OrbFerRod06}, Cai \cite{Cai07}, Brown \textit{et
al.} \cite{BroSonMcG97} an Stock and Watson
\cite{StoWat98} for applications in econometrics; Kitagawa and Gersch
\cite{KitGer85} and
Gersch and Kitagawa \cite{GerKit85} for applications in signal processing;
Hoover \textit{et al.} \cite{Hooetal98} and Ramsay and Silverman
\cite{RamSil05} for applications
in longitudinal and functional data analysis. Most of the
aforementioned literature on model \eqref{eqmodel} focused on
parameter estimation. However, it seems that the important issue of
model validation or specification of \eqref{eqmodel} have received
little attention.

For varying coefficient models of i.i.d. samples, Fan, Zhang and Zhang
\cite{FanZhaZha01} proposed the generalized likelihood ratio test
(GLRT) as a
general rule for nonparametric specification; see also Dette \cite
{Det99} for
a closely related earlier test based on nonparametric analysis of
variance (ANOVA). We also refer to the excellent review paper of Fan
and Jiang \cite{FanJia07} and the references cited therein for a more detailed
discussion of the GLRT and related tests. The GLRT has three major
advantages. First, it is of simple and intuitively appealing form. For
instance, consider testing
%
%e2 #&#
%
\begin{equation}
\label{eqsimplenull} H_0 \dvt\betav(\cdot)=\betav_0(
\cdot)\quad\longleftrightarrow\quad  H_a \dvt\betav(\cdot)\neq
\betav_0(\cdot),
\end{equation}
where $\betav_0(\cdot)$ is a known function on $[0,1]$. Then the GLRT
statistic is proportional to $(\RSS_0-\RSS_a)/\RSS_0$, where $\RSS_0$
and $\RSS_a$ are residual sum of squares under the null and alternative
hypothesis, respectively. Hence, it is similar in form to the classic
analysis of variance. Second, the GLRT is powerful to apply. Fan, Zhang
and Zhang \cite{FanZhaZha01} showed that the GLRT can detect local alternatives
with the optimal rate in the sense of Ingster \cite{Ing93}. Third, the test
is asymptotically nuisance parameter free; known as the Wilks
phenomenon. The Wilks phenomenon insures that the residual wild
bootstrap, that is, drawing i.i.d. samples from the centered empirical
distribution of the residuals, is asymptotically consistent for the
inference. In fact, the Wilks phenomenon is shown to hold for a wide
range of nonparametric models when testing under the GLRT. See, for
instance, Fan and Jiang \cite{FanJia05} for additive models and Fan
and Huang~\cite{FanHua05} for varying coefficient partially linear models. For
state-domain nonparametric regression for stationary time series, Hong
and Lee \cite{HonLee} showed that the Wilks phenomenon continue to
hold when
the errors are conditionally homogeneous.

%Investigating the behavior of GLRT for regression models under
%dependence is clearly an important problem. For instance, in
%functional and longitudinal data analysis, it is common practice to
%compare the $\RSS$ of two competing models to check goodness of fit;
%see for instance Ramsay and Silverman (2005). However, due to the
%complicated dependence structure of the data, the comparison is
%usually performed in a heuristic way. The importance of the problem
%has led to several conjectures on GLRT for dependent data. For
%instance, Chapter 9 of Fan and Yao (2003) conjectured that the results
%of Fan, Zhang and Zhang (2001) should hold for dependent data under
%some mixing conditions. Cai (2007) conjectured that the residual wild
%bootstrap is consistent for time varying coefficient model
%that the residual wild bootstrap is consistent for
%functional-coefficient regression models for nonlinear time series.

In this paper, we shall prove that the Wilks phenomenon is sensitive to
either conditional heteroscedasticity of the errors,
non-stationarity or temporal dependence in model \eqref{eqmodel}. In
particular, the Wilks phenomenon fails for model \eqref{eqmodel} even
when the errors and regressors are stationary and conditionally
homogeneous. The latter finding is drastically different from the state
domain regression case in Hong and Lee \cite{HonLee} where the Wilks
phenomenon is shown to hold when the errors are conditionally
homogeneous. As a consequence, the residual wild bootstrap fails for
model \eqref{eqmodel} under dependence since the latter bootstrap
generates (conditional) i.i.d. samples and hence mimics the Wilks type
asymptotic behavior. A new robust methodology is needed when performing
model specification for \eqref{eqmodel} under dependence and non-stationarity.

According to a result on Gaussian quadratic form approximation to the
GLRT, we shall propose in this paper a new wild bootstrap method for
the nonparametric specification of model \eqref{eqmodel}. The latter
bootstrap is shown to be consistent under non-stationarity and
dependence. We further discover that the GLRT, though fails to be
asymptotically pivotal, retains the minimax rate of local alternative
detection under
weak dependence and non-stationarity. Hence, the GLRT with the robust
wild bootstrap is powerful to apply. Note that Zhou and Wu \cite{ZhoWu10}
discussed simultaneous confidence band (SCB) construction for model
\eqref{eqmodel} which could be used for model specification. However,
the SCB can detect local alternatives with inferior rates than that of
the GLRT and hence is not a powerful tool for specification.

It is known that nonparametric specification is sensitive to the choice
of smoothing bandwidth. To alleviate the problem, Horowitz and Spokoiny
\cite{HorSpo01} and Fan, Zhang and Zhang \cite{FanZhaZha01}, among
others, proposed to
maximize the test statistic over a wide range of bandwidths. However,
for the GLRT test, the asymptotic behavior of the resulting statistic
is unknown even for i.i.d. samples, which hampers the application of
the latter test. It is worth mentioning that Zhang \cite{Zha03}
derived the
asymptotic null distribution of the maximum test for a bounded number
of bandwidths. On the other hand, M\"{u}ller \cite{Mul07} suggested to
average the GLRT over a range of bandwidths as an alternative to the
maximum test. The latter suggestion stems from surprising results, such
as Lehmann \cite{Leh06}, that the averaged likelihood ratio test can
be more
powerful than the maximum likelihood ratio test for complex
alternatives. In this paper, we shall propose to use the averaged test
for the specification of model \eqref{eqmodel} to alleviate the
sensitivity of the test to the choice of bandwidth. We derive the
asymptotic distribution and the local power of the averaged test. It is
found that the averaged test is asymptotically at least as powerful as
the best test based on a single bandwidth regardless of the shape of
the alternative, the non-stationary dependence structure of the data or
the kernel function. Our finding is potentially interesting for a wide
range of nonparametric specification problems.

Recently, there have been many results on modeling non-stationary time
series from the spectral domain. See, for instance, Dahlhaus \cite{Dah97},
Nason \textit{et al.} \cite{NasvonKro00} and Ombao \textit{et al.}
\cite{OmbvonGuo05}, among others. At the same
time, there is a great recent interest in specification of
non-stationary time series in the spectral domain. Examples include,
among others, Dahlhaus \cite{Dah09}, Neumann and von Sachs \cite{Neuvon97},
Paparoditis \cite{Pap09,Pap10}, Sergides and Paparoditis \cite
{SerPap09} and Dette \textit{et al.} \cite{DetPreVet11}. However, for the varying coefficient
regression \eqref
{eqmodel}, models from the spectral domain do not seem to be directly
useful for an asymptotic theory. In this paper, we shall adopt the time
domain modeling of locally stationary time series in Zhou and Wu
\cite{ZhoWu09}. The latter framework and the associated dependence measures
directly facilitate the theory of the current paper.

The rest of the paper is organized as follows. Section \ref{secpre}
introduces the GLRT statistic and the non-stationary time series models
for the error and regressor series. In Section \ref{secresults}, we
shall derive the asymptotic null distribution and local power of the
GLRT for parametric and semi-parametric null hypotheses. A detailed
discussion on the failure of the Wilks phenomenon is included. In
Section \ref{secnewtest}, we shall introduce the averaged test and
the corresponding robust bootstrap and investigate their asymptotic
behavior. In Section \ref{secsimu}, we shall construct a monte carlo
experiment to study the finite sample accuracy of the proposed averaged
test. Proofs of the asymptotic results are placed in Section \ref{secproofs}.

%Thanks to their flexibility and interpretability, varying coefficient
%models has gained great popularity and undergone a fast development
%over the last two decades. Applications range from longitudinal and
%functional data analysis to economics and finance. There is a huge
%statistical and econometric literature on the topic, see for instance
%Robinson (1989), Hastie and Tibshirani (1993) and the excellent review
%of Fan and Zhang (2008) and the references cited therein.

%s2 #&#
\section{Preliminaries}\label{secpre}

%s2.1 #&#
\subsection{The GLRT statistics}\label{subsecglrt}

Consider the testing problem \eqref{eqsimplenull}.
The GLRT compares the residual sum of squares (RSS) under the null
and alternative hypotheses, and a large difference indicates
violation of the null. We refer to Fan, Zhang and Zhang \cite
{FanZhaZha01} for a
detailed derivation of the statistic. Specifically, the GLRT statistic
%
%e3 #&#
%
\begin{equation}
\label{eqglrtoriginal} \lambda_n=\frac{n}{2}\log
\frac{\RSS_0}{\RSS_a}\approx-\frac{n}{2}\frac{\RSS_a-\RSS
_0}{\RSS_0},
\end{equation}
where $\RSS_0=\sum_{i=1}^n(y_i-\mathbf{x}_i^\top\betav_0(t_i))^2$ is
the RSS under the null hypothesis and $\RSS_a=\sum
_{i=1}^n(y_i-\mathbf{
x}_i^\top\hat{\betav}(t_i))^2$ is the RSS under the nonparametric
alternative. Here $\hat{\betav}(\cdot)$ is the local linear kernel
estimate of $\betav(\cdot)$ (Fan and Gijbels, \cite{FanGij96}),
which is
defined as
%
%e4 #&#
%
\begin{eqnarray}
\label{eqquantile} \bigl(\hat{\betav}_{b_n}(t),\hat{
\betav}'_{b_n}(t)\bigr) = \mathop{\argmin}_{\eta_0, \eta_1 \in\R^p} \sum
_{i=1}^n \bigl(y_i-
\mathbf{x}^{\top}_i\eta_0 -\mathbf{x}_i^{\top}
\eta_1(t_i-t) \bigr)^2 K_{b_n}(t_i-t),
\end{eqnarray}
where $K$ is a kernel function, $b_n > 0$ is the bandwidth, and
$K_{c}(\cdot) = K(\cdot/c)$, $c > 0$. Throughout this paper, we shall
always assume that the kernel $K \in\mathcal{ K}$, the collection of
symmetric density functions $K$ with support $[-1, 1]$ and $K \in
\mathcal{ C}^1 [-1,1]$. Define
\[
\mathbf{S}_{n, l}(t) = (nb_n)^{-1}\sum
_{i=1}^n \mathbf{x}_i\mathbf{
x}^{\top}_i \bigl[(t_i-t)/b_n
\bigr]^l K_{b_n}(t_i-t)
\]
for $l = 0, 1, \ldots,$ where $0^0:=1$, and
\[
\mathbf{R}_{n,l}(t)=(nb_n)^{-1}\sum
_{i=1}^n\mathbf{x}_i y_i
\bigl[(t_i-t)/b_n\bigr]^lK_{b_n}(t_i-t).
\]
Let $\hat{\etav}_{b_n}(t) = (\hat{\betav}^{\top}_{b_n}(t),
b_n(\hat{\betav}'_{b_n}(t))^{\top} )^{\top}$. Then it can be shown
that (Fan and Gijbels, \cite{FanGij96})
%
%e5 #&#
%
\begin{equation}
\label{eqsol} \hat{\etav}_{b_n}(t)=\pmatrix{ \mathbf{S}_{n,0}(t) & \mathbf{S}_{n,1}^{\top}(t)
\vspace*{2pt}\cr
\mathbf{S}_{n,1}(t) & \mathbf{S}_{n,2}(t)
}^{-1}\pmatrix{
\mathbf{R}_{n,0}(t)
\vspace*{2pt}\cr
\mathbf{R}_{n,1}(t)} :=
\mathbf{S}^{-1}_{n}(t)\mathbf{R}_{n}(t).
\end{equation}
We shall omit the subscript $b_n$ in $\hat{\etav}$, $\hat{\betav}$
and $\hat{\betav}'$ hereafter if no confusion will be caused.

% =2\hat{\betav}_{b_n/\sqrt{2}}(t)-\hat{\betav}_{b_n}(t).
%The bias of $\tilde\betav(t)$ is of order $o(b_n^2)$ and is
%asymptotically negligible. Implementing \eqref{eqjackknife} is
%asymptotically equivalent to using the 4th order kernel $K^*(x) = 2

%s2.2 #&#
\subsection{Locally stationary time series models}
Throughout this paper, we shall assume that both $(\mathbf{x}_i)$ and
$(\varepsilon_i)$ belong
to a general class of locally stationary time series in the sense of
Zhou and Wu \cite{ZhoWu09} as follows,
%
%e6 #&#
%
\begin{eqnarray}
\label{eqnonstatioanrymodelce} \mathbf{x}_i&=&\mathbf{G}
\bigl(t_i,(\ldots, \epsilon_{i-1}, \epsilon_i)
\bigr),\qquad  i=1,2,\ldots,n,
\nonumber
\\[-8pt]
\\[-8pt]
\nonumber
\varepsilon_i&=&H\bigl(t_i, (\ldots,
\xi_{i-1}, \xi_i)\bigr) V\bigl(t_i, (\ldots,
\epsilon_{i-1}, \epsilon_i)\bigr),\qquad i=1,2,\ldots,n,
\end{eqnarray}
where $\mathbf{G}(\cdot)=(G_1,G_2,\ldots,G_p)^{\top}(\cdot)$,
$(\epsilon_i)_{i \in\Z}$ are i.i.d., $(\xi_i)_{i \in\Z}$ are also i.i.d.
and $(\epsilon_i)_{i \in\Z}$ is independent of $(\xi_i)_{i \in
\Z}$. Let $\FF_i=(\ldots, \epsilon_{i-1}, \epsilon_i)$ and
$\GG_i=(\ldots, \xi_{i-1}, \xi_i)$. We assume that
\[
\E\bigl(H(t,\GG_i)\bigr)=0 \quad\mbox{and}\quad \var\bigl(H(t,
\GG_i)\bigr)=1,
\]
almost surely for all $t\in[0,1]$, in which case $V^2(t_i,\FF_i)$ is the
conditional variance of $\varepsilon_i$ given~$\FF_i$.

It is clear from \eqref{eqnonstatioanrymodelce} that $(\mathbf{
x}_i)$ and $(\varepsilon_i)$ are non-stationary. Formulation
\eqref{eqnonstatioanrymodelce} can be interpreted as physical
systems with $\mathcal{ F}_i$ and $\GG_i$ being the inputs and
$\mathbf{
x}_i$, $\varepsilon_i$ being the outputs, respectively, and $\mathbf{
G}$, $H$ and $V$ being the transforms or filters that represent the
underlying physical mechanism. By allowing~$\mathbf{G}$, $H$ and $V$
varying smoothly with respect to $t$, we have local stationarity of
$(\mathbf{x}_i)$ and $(\varepsilon_i)$. See also Zhou and Wu \cite
{ZhoWu09} for
more discussions. The above formulation of covariates and error
processes is very general and includes many settings in the existing
time series regression literature as special cases. To help
understand the formulation, we shall consider the following three
cases:
\begin{longlist}[(a)]
\item[(a)] (I.i.d. model). Assume that $\mathbf{x}_i=\mathbf
{G}_0(\epsilon_i)$
and $\varepsilon_i=H_0(\xi_i)$. Then $(\mathbf{x}^{\top}_i,
\varepsilon_i)_{i=1}^n$ is a random sample and
$(\varepsilon_i)_{i=1}^n$ is independent of $(\mathbf{x}_i)_{i=1}^n$.
This type of design was discussed extensively in Fan, Zhang and
Zhang \cite{FanZhaZha01} and Fan and Jiang \cite{FanJia07}, among others.

\item[(b)](Exogenous model). In \eqref{eqnonstatioanrymodelce}, we
assume that $V(t_i, \FF_i)=V_0(t_i)$. In this case, the
regressors and errors are two independent locally stationary
processes. Under further restrictions on the processes, this type of
model was studied in Robinson \cite{Rob89}, Orbe \textit{et al.}
\cite{OrbFerRod05,OrbFerRod06} among
others.

\item[(c)](Endogenous model). Assume \eqref{eqnonstatioanrymodelce}.
Note that in this case the errors are correlated with the regressors
since they both depend on inputs $\FF_i$. This type of model is
suitable when the errors exhibit heteroscedasticity with respect to
time and the regressors. When $\mathbf{x}_i$ and $H(t,\GG_i)$ are
stationary, the case was considered in Cai \cite{Cai07} among others.
\end{longlist}

Write $\chi_i=(\epsilon_i,\xi_i)^{\top}$ and
$\RR_i=(\ldots,\chi_{i-1},\chi_i)$. For a generic locally stationary
time series $Z_i=\mathbf{L}(t_i,\RR_i)$. The strength of the temporal
dependence in $\{Z_i\}$ can be measured by how strongly the `current'
observation of the time series, $Z_i$, is influenced by the innovation
$\chi_0$ which occurred $i$ steps ahead. More specifically, we can define
%
%e7 #&#
%
\begin{eqnarray}
\label{eqdep} \delta_{p}(\mathbf{L},k)=\sup_{0\le t\le1}\bigl\|
\mathbf{L}(t,\RR_k)-\mathbf{L}\bigl(t,\RR^*_{k}\bigr)
\bigr\|_p \qquad\mbox{where } \RR_k^*=\bigl(\RR_{-1},
\chi^*_0,\chi_1,\chi_2,\ldots,
\chi_i\bigr)
\end{eqnarray}
and $\{\chi_i^*\}$ is an i.i.d. copy of $\{\chi_i\}$. Implementing the
idea of coupling, $\delta_p(\mathbf{L},k)$ measures the effect of
$\chi_0$
in generating observations that are $k$ steps away. Therefore, if
$\delta_p(\mathbf{L},k)$ decays fast as $k$ gets large, short range
dependence is implied. We refer to Zhou and Wu \cite{ZhoWu09} for more
discussions and examples on the above dependence measures.

%s3 #&#
\section{Asymptotic results}\label{secresults}
For a family of stochastic processes $(\mathbf{L}(t, \RR_i))_{i
\in\Z}$, we say that it is $\mathcal{ L}^q$ stochastic Lipschitz
continuous on $[0,1]$ if $\sup_{0 \le s < t\le1} [\|\mathbf
{L}(t,\mathcal{
R}_0) - \mathbf{L}(s,\mathcal{ R}_0)\|_q / (t-s)] < \infty$. Denote by
$\mathrm{Lip}_q$ the collection of such systems. Let $\mathcal{ U}^{p}$ be the
collection of processes $(\mathbf{L}(t, \RR_i))_{i
\in\Z}$ such that $\|\mathbf{L}(t, \RR_0))\|_p<\infty$ for all
$t\in
[0,1]$. Let $\mathcal{ C}^l \mathcal{ I}$,
$l \in\N$, be the collection of functions that have $l$th order
continuous derivatives on the interval $\mathcal{ I} \subset\R$. We
shall make the following assumptions:

\begin{longlist}[(A1)]
\item[(A1)]
Let $M(t)$ be the $p\times p$ matrix with $(i,j)$th entry $m_{ij}(t)
= \E[G_i(t,\FF_0) G_j(t,\FF_0)]$. Assume that the smallest
eigenvalue of $M(t)$ is bounded away from $0$ on $[0,1]$ and
$M(t)\in\mathcal{ C}^2[0,1]$.

\item[(A2)]
$\mathbf{G}(t,\FF_i)\in\mathcal{ U}^{32} \cap \mathrm{Lip}_2$ for some $r>0$.

\item[(A3)] $\mathbf{U}(t,\RR_i) := \mathbf{G}(t, \FF_i)V(t,\FF
_i) H(t, \GG_i) \in
\mathcal{ U}^4 \cap \mathrm{Lip}_2$.

\item[(A4)]
$\sum_{k=0}^{\infty}\delta_{32}(\mathbf{G},k)<\infty$.

\item[(A5)]
$\delta_{4}(V,k)+\delta_{4}(H,k)=\mathrm{O}((k+1)^{-2})$.

\item[(A6)]
$\delta_4(\mathbf{U},k)=\mathrm{O}(\chi^k)$ for some $\chi\in(0,1)$.

\item[(A7)]
The smallest eigenvalue of $\Lambda(t)$ is bounded away from $0$ on
$[0,1]$, where
%
%e8 #&#
%
\begin{equation}
\label{eqlrv} \Lambda(t)=\sum_{i=-\infty}^{\infty}
\operatorname{cov} \bigl(\mathbf{U}(t,\RR_0),\mathbf{U}(t,
\RR_i) \bigr).
\end{equation}

\item[(A8)]
The coefficient functions $\beta_j(\cdot) \in\mathcal{ C}^2[0,1]$,
$j =
1, \ldots, p$.
\end{longlist}

A few remarks on the regularity conditions are in order. Conditions
(A1), (A2) and (A4) insures local stationarity and short memory of the
regressor process $\mathbf{x}_i$. The existence of the $32$rd moment is
for technical convenience only and may be relaxed. The eigenvalue
constraint in condition (A1) insures the non-singularity of the design.
Conditions (A3), (A5) and (A6) guarantees the smoothness and short
range dependence of the error process $\varepsilon_i$. Furthermore,
condition (A7) means that the asymptotic covariance matrix of $\hat
{\betav}(t)$ is non-singular.

%s3.1 #&#
\subsection{The null distributions}\label{secnull}
%
%th1 #&#
%
\begin{theorem}\label{thm1}
Assume that condition \textup{(A)} holds and that $nb_n^{9/2}=\mathrm{O}(1)$ and
$nb_n^4/(\log n)^6\rightarrow\infty$. Then under $H_0$, we have
%
%e9 #&#
%
\begin{eqnarray}
\label{eqnull} \hspace*{-15pt}\sqrt{b_n}\biggl\{2\lambda_n+
\frac{\tilde{K}(0)}{b_n\mathcal{ V}}\int_{0}^1\tr\bigl[H(t)\bigr]
\,\mathrm{d}t+\frac{nb_n^4\mu_2^{2}}{4\mathcal{ V}}\int_0^1\bigl[
\betav''(t)\bigr]^\top M(t)
\betav''(t) \,\mathrm{d}t\biggr\}\Rightarrow N\bigl(0,
\sigma^2/\mathcal{ V}^2\bigr),
\end{eqnarray}
where
\[
\sigma^2=\int_{\R}\tilde{K}^2(t) \,
\mathrm{d}t\int_{0}^1\tr\bigl[H(t)^2
\bigr] \,\mathrm{d}t,
\]
$\tilde{K}(\cdot)=K\ast K(\cdot)-2K(\cdot)$, $H(\cdot)=\Lambda
^{1/2}(\cdot)M^{-1}(\cdot)\Lambda^{1/2}(\cdot)$, $\mathcal{
V}=\int_0^1\E[
V(t, \FF_0)]^2 \,\mathrm{d}t$, $\mu_2=\int_{-1}^1x^2K(x) \,\mathrm
{d}x$, `$\ast$' is the
convolution operator and `$\tr$' denotes the trace of a matrix.
\end{theorem}

Theorem \ref{thm1} reveals the asymptotic behavior of the GLRT for a
very wide class of predictor and error processes. In particular, the
latter Theorem explains when and why the Wilks phenomenon fails. In the
following, we will consider four special cases to see how endogeneity,
non-stationarity and temporal dependence influence the Wilks
phenomenon. To simplify the discussion, we will assume in the examples
below that the asymptotic bias effect, $\frac{nb_n^4\mu
_2^{2}}{4\mathcal{
V}}\int_0^1[\betav''(t)]^\top M(t)\betav''(t) \,\mathrm{d}t$, is
asymptotically
negligible in \eqref{eqnull}. In practice, the latter task can be
achieved by pre-whitening. We will discuss bias reduction techniques
for GLRT in Section~\ref{secbandwidth}.

%ex1 #&#
%
\begin{example}[(I.i.d. sample without endogeneity)]\label{ex1}
Consider the case when $\mathbf{x}_i=\mathbf{G}(\epsilon_i)$ and
$\varepsilon_i=CH(\zeta_i)$, where $C$ is a positive constant. In
this case, the
covarites and errors are two independent i.i.d. sequences and the
conditions in Fan, Zhang and Zhang \cite{FanZhaZha01} are satisfied.
Note that
$\mathcal{ V}=C^2$, $\Lambda(t)=M(t)C^2$ and $H(t)=C^2\mathbf{I}_p$, where
$\mathbf{I}_p$ is the $p\times p$ identity matrix. In particular,
%
%e10 #&#
%
\begin{equation}
\label{eqpivotal} \int_{0}^1\tr\bigl[H(t)\bigr]
\,\mathrm{d}t/\mathcal{ V}=p\quad \mbox{and}\quad \int_{0}^1
\tr\bigl[H(t)^2\bigr] \,\mathrm{d}t/\mathcal{ V}^2=p
\end{equation}
in \eqref{eqnull}. Hence, it is easy to check that
\[
\sqrt{b_n}\biggl\{2\lambda_n+\frac{p\tilde{K}(0)}{b_n}\biggr\}
\Rightarrow N\biggl(0,p\int_{\R}\tilde{K}^2(t) \,
\mathrm{d}t\biggr),
\]
which coincides with Theorem 5 of Fan, Zhang and Zhang \cite
{FanZhaZha01} and the
Wilks phenomenon holds.
\end{example}

%ex2 #&#
%
\begin{example}[(The effect of temporal dependence)]\label{ex2}
In this case $\mathbf{x}_i=\mathbf{G}(\FF_i)$ and $\varepsilon
_i=CH(\GG_i)$,
where $C$ is a positive constant. Hence, $\{\mathbf{x}_i\}$ and $\{
\varepsilon_i\}$ are two stationary processes which are independent of
each other. In particular, neither endogeneity nor non-stationary is
assumed in the model. It is easy to see that, in this case,
%
%e11 #&#
%
\begin{equation}
\label{eqstationarylrt} \Lambda(t)=C^2\sum
_{i=-\infty}^{\infty}\E\bigl[\mathbf{G}(\FF_0)
\mathbf{G}^{\top
}(\FF_i)\bigr]\E\bigl[H(\GG_0) H(
\GG_i)\bigr],
\end{equation}
$\mathcal{ V}=C^2$ and $M(t)=\E[\mathbf{G}(\FF_0) \mathbf{G}^{\top
}(\FF_0)]$. An
important observation is that
\begin{eqnarray*}
&&\!\!\int_{0}^1\!\tr\bigl[H(t)\bigr] \,\mathrm{d}t/
\mathcal{ V}
=\tr\Biggl(\bigl\{\E\bigl[\mathbf{G}(\FF_0)
\mathbf{
G}^{\top}(\FF_0)\bigr]\bigr\}^{-1}\!\sum
_{i=-\infty}^{\infty}\!\E\bigl[\mathbf{G}(\FF_0)
\mathbf{ G}^{\top}(\FF_i)\bigr]\E\bigl[H(\GG_0)
H(\GG_i)\bigr] \Biggr),
\\
&&\!\!\int_{0}^1\!
\tr\bigl[H(t)^2\bigr] \,\mathrm{d}t/\mathcal{ V}^2
=\tr
\Biggl(\Biggl[\!\bigl\{\E\bigl[\mathbf{G}(\FF_0) \mathbf{G}^{\top}(
\FF_0)\bigr]\bigr\}^{-1}\!\sum_{i=-\infty}^{\infty}\!
\E\bigl[\mathbf{G}(\FF_0) \mathbf{G}^{\top}(
\FF_i)\bigr]\E\bigl[H(\GG_0) H(\GG_i)\bigr]
\Biggr]^{\!2} \Biggr)
\end{eqnarray*}
are no longer nuisance parameter free compared with the results in
\eqref{eqpivotal}. As a consequence, the Wilks phenomenon fails to
hold in this case. Additionally, it is easy to see that the latter loss
of pivotality is due to the fact that the summands in \eqref
{eqstationarylrt} are generally nonzero for $i\neq0$, which is
caused by the temporal dependence. Indeed, if the summands are zero for
$i\neq0$ in \eqref{eqstationarylrt}, then $\Lambda(t)=C^2\E[\mathbf{
G}(\FF_0) \mathbf{G}^{\top}(\FF_0)]$ and we have \eqref
{eqpivotal}. Like
in many pivotal tests such as the Wald test, the term $\RSS_0/n\approx
\mathcal{ V}$ in the GLRT serves as a scaling device which cancels out the
variance factor in $\RSS_1-\RSS_0$ and makes the test pivotal in the
i.i.d. case. However, as shown above, $\RSS_0/n$ fails to fulfill the
latter scaling task under dependence.
\end{example}

%ex3 #&#
%
\begin{example}[(The effect of non-stationarity)]\label{ex3}
Let $\mathbf{x}_i=\mathbf{G}(t_i,\epsilon_i)$ and $\varepsilon
_i=V(t_i)H(t_i,\zeta_i)$. Here $\{\mathbf{x}_i\}$ and $\{\varepsilon
_i\}$
are two independent but non-stationary sequences which are independent
of each other. In this case, we have
%
%e12 #&#
%
\begin{eqnarray}
\label{eqnon-stati} \int_{0}^1\tr\bigl[H(t)
\bigr] \,\mathrm{d}t/\mathcal{ V}=p\quad \mbox{and}\quad \int_{0}^1
\tr\bigl[H(t)^2\bigr] \,\mathrm{d}t/\mathcal{ V}^2=p
\frac{\int_0^1V^4(t) \,\mathrm{d}t}{(\int_0^1V^2(t) \,\mathrm{d}t)^2}.
\end{eqnarray}
Note that the second term in \eqref{eqnon-stati} depends on the
time-varying variance $V^2(t)$ and hence the Wilks phenomenon fails to
hold in this case. Additionally, observe that $\frac{\int_0^1V^4(t)
\,\mathrm{d}t}{(\int_0^1V^2(t) \,\mathrm{d}t)^2}\ge1$ and the
equation holds if and only if
$V(t)$ is a constant function. Compared with the results in \eqref
{eqpivotal}, we conclude that, in this case, non-stationarity in the
errors tends to inflate the variance of GLRT. Furthermore, if $\{
\varepsilon_i\}$ has constant variance, then the Wilks phenomenon holds
even if $\{\mathbf{x}_i\}$ is a non-stationary sequence.
\end{example}

%ex4 #&#
%
\begin{example}[(The effect of endogeneity)]\label{ex4}
Suppose that $\mathbf{x}_i=\mathbf{G}(\epsilon_i)$ and $\varepsilon
_i=V(\epsilon_i)H(\zeta_i)$. In this case $\{\mathbf{x}_i\}$ and $\{
\varepsilon_i\}$ are two i.i.d. sequences which are dependent of each
other. We obtain
\begin{eqnarray*}
\int_{0}^1\tr\bigl[H(t)\bigr] \,\mathrm{d}t/
\mathcal{ V}&=&\tr\bigl(\bigl\{\E\bigl[\mathbf{G}(\epsilon
_0)\mathbf{
G}^{\top}(\epsilon_0)\bigr]\bigr\}^{-1}\E\bigl[
\mathbf{G}(\epsilon_0)\mathbf{G}^{\top
}(\epsilon_0)V^2(
\epsilon_0)\bigr] \bigr)/\E\bigl[V^2(\epsilon_0)
\bigr],
\\
\int_{0}^1\tr\bigl[H(t)^2
\bigr] \,\mathrm{d}t/\mathcal{ V}^2&=&\tr\bigl(\bigl[\bigl\{\E
\bigl[
\mathbf{G}(\epsilon_0)\mathbf{G}^{\top}(\epsilon_0)
\bigr]\bigr\}^{-1}\E\bigl[\mathbf{G}(\epsilon_0)
\mathbf{G}^{\top
}(\epsilon_0)V^2(
\epsilon_0)\bigr]\bigr]^2 \bigr)/\bigl(\E
\bigl[V^2(\epsilon_0)\bigr]\bigr)^2.
\end{eqnarray*}
Note that if $\E[\mathbf{G}(\epsilon_0)\mathbf{G}^{\top}(\epsilon
_0)V^2(\epsilon_0)]=\E[\mathbf{G}(\epsilon_0)\mathbf{G}^{\top
}(\epsilon_0)]\E
[V^2(\epsilon_0)]$, then we have \eqref{eqpivotal} and hence the Wilks
phenomenon. Due to the dependence of $\mathbf{G}(\epsilon_0)$ and
$V(\epsilon_0)$, the latter factorization generally fails and hence the
Wilks phenomenon fails to hold in this case.
\end{example}

In many real applications, one is interested in specifying a component
of $\betav(\cdot)$. For instance, one may want to test whether $\beta
_j(\cdot)$ is significantly different from zero. This leads us
to\vadjust{\goodbreak}
consider the following hypothesis testing problem where both $H_{01}$
and $H_{a1}$ are nonparametric:
%
%e13 #&#
%
\begin{eqnarray}
\label{eqseminull} H_{01} \dvt\betav^{(1)}(\cdot)=
\betav^{(1)}_{0}(\cdot)\quad \longleftrightarrow \quad H_{a1}
\dvt\betav^{(1)}(\cdot)\neq\betav^{(1)}_{0}(\cdot),
\end{eqnarray}
where
\begin{eqnarray*}
\betav(t)=\pmatrix{\betav^{(1)}(t)
\vspace*{2pt}\cr
\betav^{(2)}(t)},\qquad
\betav_0(t)=\pmatrix{
\betav^{(1)}_{0}(t)
\vspace*{2pt}\cr
\betav^{(2)}_{0}(t)
}
\quad\mbox{and}\quad \mathbf{x}_i=\pmatrix{
\mathbf{x}^{(1)}_{i}
\vspace*{2pt}\cr
\mathbf{x}^{(2)}_i
}
,
\end{eqnarray*}
$\betav^{(1)}(t)$, $\betav^{(1)}_{0}(t)$ and $\mathbf{x}_{i}^{(1)}$ are
$p_1<p$ dimensional and $\betav^{(1)}_{0}(t)$ is a known function.
Define $y_i^*=y_i-[\betav^{(1)}_{0}(t_i)]^{\top}\mathbf{x}_{i}^{(1)}$.
Then under $H_{01}$ the functions $\beta_j(\cdot)$, $j=p_1+1,\ldots,p$
can be estimated by the local linear regression of $y^*_i$ on $\mathbf{
x}_{i}^{(2)}$ with bandwidth $b_n$. Throughout the paper we assume that
the bandwidth $b_n$ used under $H_{01}$ is the same as that under
$H_{a1}$. Asymptotic results can be easily obtained using the arguments
of the paper when the two bandwidths are different. However, the
resulting asymptotic bias and variance are much more complicated. For
the sake of presentational clarity, we will only consider the case of
equal bandwidth.

The GLRT statistic for testing $H_{01}$ against $H_{a1}$ is defined as
%
%e14 #&#
%
\begin{eqnarray}
\label{eqglrtsemi} \lambda_{1n}=\frac{n}{2}\log
\frac{\RSS_1}{\RSS_a}=\frac{n}{2} \biggl[\log\frac{\RSS_1}{\RSS
_0}-\log
\frac{\RSS_a}{\RSS_0} \biggr]\approx-\frac{n}{2}\frac{\RSS
_a-\RSS_1}{\RSS_0},
\end{eqnarray}
where $\RSS_1$ is the RSS under $H_{01}$.

Write
\begin{eqnarray*}
M(t)=\pmatrix{M_{11}(t)&M_{12}(t)
\vspace*{2pt}\cr
M_{21}(t)&M_{22}(t)
}
\quad\mbox{and}\quad \Lambda(t)=\pmatrix{
\Lambda_{11}(t)&\Lambda_{12}(t)
\vspace*{2pt}\cr
\Lambda_{21}(t)&\Lambda_{22}(t)
}
,
\end{eqnarray*}
where $M_{11}(t)$ and $\Lambda_{11}(t)$ are of dimension $p_1\times
p_1$.

Define $p\times p$ matrix $H_2(t)=\operatorname{diag}(\mathbf
{0}_{p_1},\Lambda^{1/2}_{22}(t)M^{-1}_{22}(t)\Lambda^{1/2}_{22}(t))$.
We have the
following theorem.
%
%th2 #&#
%
\begin{theorem}\label{thm2}
Assume that condition \textup{(A)} holds and that $nb_n^{9/2}=\mathrm{O}(1)$ and
$nb_n^4/(\log n)^6\rightarrow\infty$. Then under $H_{01}$, we have
\begin{eqnarray*}
\sqrt{b_n}\biggl\{2\lambda_{1n}+\frac{\tilde{K}(0)}{b_n\mathcal{
V}}\int
_{0}^1\tr\bigl[H^*(t)\bigr] \,\mathrm{d}t+
\frac{nb_n^4\mu_2^{2}}{4\mathcal{ V}}\int_0^1\Upsilon(t) \,
\mathrm{d}t\biggr\} \Rightarrow N\bigl(0,\sigma_1^2/
\mathcal{ V}^2\bigr),
\end{eqnarray*}
where $H^*(\cdot)=H(\cdot)-H_2(\cdot)$, $\Upsilon(t)=[\betav
''(t)]^\top
M(t)\betav''(t)-\{[\betav^{(2)}(t)]''\}^\top M_{22}(t)[\betav
^{(2)}(t)]''$ and
\[
\sigma_1^2=\int_{R}
\tilde{K}^2(t) \,\mathrm{d}t\int_{0}^1
\tr\bigl[\bigl\{H^*(t)\bigr\}^2\bigr] \,\mathrm{d}t.
\]
\end{theorem}

Theorem \ref{thm2} unveils the asymptotic null distribution of the
test under $H_{01}$. Following very similar arguments as those in
Examples~\ref{ex1}--\ref{ex4}, the Wilks phenomenon can be shown to be sensitive to
non-stationary, temporal dependence and endogeneity in this case as well.

Practitioners and researchers often encounter testing problems where
the null is specified up to a parametric part. For instance, one may
want to test whether $\betav(\cdot)$ is really time varying in model~\eqref{eqmodel}, which amounts to testing $\betav(\cdot)=C$ for some
unspecified constant vector $C$. Heuristically, since the convergence
rate of the local linear estimates is always slower than the $\sqrt{n}$
parametric rate, it is expected that the null distribution will not be
altered as long as we plug in a $\sqrt{n}$ consistent estimate of the
unspecified parametric part. The following discussion rigorously
confirms the intuition. Consider testing
\[
\tilde{H}_{01} \dvt\betav^{(1)}(\cdot)=\betav^{(1)}_{0}(
\cdot,\theta_0)\qquad\mbox{for some unknown } \theta_0\in\Omega
\subset\R^q,
\]
where $\{\betav^{(1)}_{0}(\cdot,\theta)\dvt\theta\in\Omega\}$ is a
parametric family of smooth functions.
Let $\tilde{y}_i^*=y_i-(\betav^{(1)}_{0})^{\top}\times (t_i,\hat{\theta
})\mathbf{
x}_{i}^{(1)}$ and $\widetilde{\RSS}_1$ be the residual sum of squares
of the local linear regression of $\tilde{y}^*_i$ on $\mathbf{
x}_{i}^{(2)}$ with bandwidth $b_n$. We shall make the following
assumptions on the parametric family $\betav^{(1)}_0(\cdot,\theta)$ and
the estimate $\hat{\theta}$:

\begin{longlist}[(B1)]
\item[(B1)] For each $t\in[0,1]$, $\betav^{(1)}_{0}(t,\theta)$ is
$\mathcal{ C}^2$ in $\theta$ in a neighborhood $\Theta$ of $\theta
_0$. Additionally,
\[
\sup_{t\in[0,1],\theta\in\Theta} \biggl\{ \biggl|\frac{\partial\betav
^{(1)}_{0}(t,\theta)}{\partial\theta} \biggr|+ \biggl|\frac{\partial^2 \betav
^{(1)}_{0}(t,\theta)}{\partial\theta^2} \biggr|
\biggr\}<\infty.
\]

\item[(B2)] Under $\tilde{H}_{01}$, $\|\hat{\theta}-\theta_0\|
_4=\mathrm{O}(1/\sqrt{n})$.
\end{longlist}

%pr1 #&#
%
\begin{proposition}\label{propequasemi}
Under $\tilde{H}_{01}$, condition \textup{(B)} and the assumptions of Theorem
\ref{thm2}, we have
%
%e15 #&#
%
\begin{equation}
\label{eqequasemi} \widetilde{\RSS}_1-\RSS_1-\mathrm{O}_\p
\bigl(\sqrt{n}b_n^2\bigr)=\mathrm{O}_\p(1).
\end{equation}
\end{proposition}
The $\mathrm{O}_\p(\sqrt{n}b_n^2)$ term on the left-hand side of \eqref
{eqequasemi} corresponds to the extra bias introduced by the
estimation error of $\theta$. And the $\mathrm{O}_\p(1)$ term on the right-hand
side of \eqref{eqequasemi} corresponds to the extra variance caused
by the latter error. Both terms are asymptotically negligible compared
to the $\mathrm{O}_\p(nb_n^4)$ bias and $\mathrm{O}_\p(1/b_n)$ variance of $\RSS_1$.
As a
consequence, the results of Theorems \ref{thm1} and \ref{thm2}
continues to hold if $\theta$ is replaced by $\hat{\theta}$.

%s3.2 #&#
\subsection{Local power of the GLRT}
%
%pr2 #&#
%
\begin{proposition}\label{proplocalpower}
Assume the alternative $H_{a,n}\dvt \betav(\cdot)=\betav_0(\cdot
)+n^{-4/9}\mathbf{f}_n(\cdot)$, where $\mathbf{f}_n(\cdot)\in
\mathcal{
C}^2[0,1]$. Further assume that $b_n=cn^{-2/9}$ for some $c>0$, that
$\int_0^1|\mathbf{f}_n''(t)| \,\mathrm{d}t=\mathrm{o}(n^{4/9})$ and that
\begin{eqnarray*}
\int_0^1\mathbf{f}_n^{\top}(t)M(t)
\mathbf{f}_n(t) \,\mathrm{d}t\rightarrow F_1,\qquad
n^{-8/9}\int_0^1\bigl[
\mathbf{f}_n''(t)\bigr]^\top
M(t)\mathbf{f}_n''(t) \,\mathrm{d}t
\rightarrow F_2
\end{eqnarray*}
for some finite constants $F_1$ and $F_2$. Then under condition \textup{(A)}, we have
\begin{eqnarray*}
\hspace*{-4pt}&&\sqrt{b_n} \biggl\{2\lambda_n+\frac{\tilde{K}(0)}{b_n\mathcal{
V}}\int
_{0}^1\tr\bigl[H(t)\bigr] \,\mathrm{d}t \biggr
\}+\frac{c^{9/2}\mu_2^{2}}{4\mathcal{ V}}\int_0^1\bigl[
\betav''(t)\bigr]^\top M(t)
\betav''(t) \,\mathrm{d}t+\frac{c^{9/2}\mu_2^{2}}{4\mathcal{ V}}F_2
-\frac{c^{1/2}}{\mathcal{ V}}F_1\\
\hspace*{-4pt}&&\quad\Rightarrow N\bigl(0,\sigma^2/
\mathcal{ V}^2\bigr).
\end{eqnarray*}
\end{proposition}

%that $\int_0^1[|\mathbf{f}_n(t)|^2 \,\mathrm{d}t+n^{-8/9}\int_0^1|
% \,\mathrm{d}t<C$ for some finite constant $C$

When the errors and regressors are weakly dependent locally stationary
time series, Proposition \ref{proplocalpower} claims that the GLRT
can still detect local alternatives with the optimal rate $\mathrm{O}(n^{-4/9})$
in the sense of Ingster \cite{Ing93}. As a consequence, the GLRT is powerful
to apply for nonparametric model validation of model \eqref{eqmodel}
under non-stationarity and dependence. However, it should be noted that
the GLRT may not be the most powerful among all rate optimal tests. In
the literature, among other examples, Zhang and Dette \cite{ZhaDet04}
discovered
that other tests may yield smaller variance than the GLRT for
independent samples. From Proposition \ref{proplocalpower}, the
asymptotic local power of the GLRT with level $\alpha$
%
%e16 #&#
%
\begin{eqnarray}
\label{powerglrt} \beta_\alpha(c)=\Phi(R_1-z_{1-\alpha})\qquad
\mbox{where } R_1=\frac
{c^{1/2}F_1-c^{9/2}\mu_2^2F_2/4}{\sigma},
\end{eqnarray}
$\Phi(\cdot)$ and $z_{1-\alpha}$ denote the cumulative distribution
function and the $1-\alpha$ quantile of the standard normal
distribution. Assume that $F_1\neq0$ and $F_2\neq0$, then simple
calculations show that the bandwidth which maximizes the above power is
\[
\tilde{b}_n=\tilde{c}n^{-2/9}\qquad \mbox{where } \tilde{c}=
\biggl(\frac
{4F_1}{9\mu_2^2F_2} \biggr)^{1/4}.
\]

%re1 #&#
%
\begin{remark}
A typical example which satisfies $F_1\neq0$ and $F_2\neq0$ is when
$\mathbf{f}_n(t)=a_n\mathbf{f}(a_n^2(t-t_0))$, where $\mathbf{f}\in
\mathcal{
C}^2[-1,1]$, $t_0\in(0,1)$ and $a_n=n^{1/9}$. Simple calculations show that
%
%e17 #&#
%
\begin{eqnarray}
\label{eqinte} F_1=\int_{-1}^1
\mathbf{f}^{\top}(t)M(t_0)\mathbf{f}(t) \,\mathrm{d}t,\qquad
F_2=\int_{-1}^1\bigl[
\mathbf{f}''(t)\bigr]^\top
M(t_0)\mathbf{f}''(t) \,\mathrm{d}t.
\end{eqnarray}
Hence $F_1\neq0$ and $F_2\neq0$ as long as the corresponding terms in
\eqref{eqinte} are nonzero.
\end{remark}

%s4 #&#
\section{Tests for locally stationary time series}\label{secnewtest}

%s4.1 #&#
\subsection{The test}

Consider the testing problem \eqref{eqsimplenull}. Two important
observations lead to the following modifications of the original GLRT
when testing for non-stationary time series. First, as shown in
Examples~\ref{ex2}--\ref{ex4}, the denominator $\RSS_0/n$ is redundant when testing for
non-stationary time series. Second, as we discussed in the
\hyperref[secintr]{Introduction}, averaging the test over a range of
bandwidths can reduce
the sensitivity of the test with respect to the selection of bandwidth
and may also gain power over tests based on a single (optimal)
bandwidth. Based on the above discussions, we suggest using the
following averaged test when specifying model \eqref{eqmodel} for
non-stationary time series:
%
%e18 #&#
%
\begin{equation}
\label{eqtest} \lambda^*_n=\int_{c_{\min}}^{c_{\max}}
\bigl(\RSS_{0}-\RSS_{a}\bigl(zn^{-\gamma
}\bigr)\bigr)
\,\mathrm{d}z,
\end{equation}
where $\RSS_{a}(b)$ is the RSS under $H_a$ when bandwidth is chosen as
$b$, $0<c_{\min}<c_{\max}<\infty$. Large $\lambda^*_n$ indicates
evidence against $H_0$. In the literature, nonparametric ANOVA tests
ignoring the denominator were first proposed in Dette \cite{Det99} for
independent samples. Dette and Hetzler~\cite{DetHet07} also considered averaged
nonparametric specification tests over a range of bandwidths. The
following theorem derives the asymptotic null distribution of the
averaged test.

%th3 #&#
%
\begin{theorem}\label{thmnewtest}
Assume that condition \textup{(A)} holds and that $2/9\le\gamma<1/4$. Then
under $H_0$, we have %and that %$w(\cdot)\in\mathcal{ C}^1[c_{\min},c_{
\begin{eqnarray*}
&&\sqrt{n^{-\gamma}} \biggl\{\lambda^*_n+n^{\gamma}
\tilde{K}(0)\bigl[\log(c_{\max})-\log(c_{\min})\bigr]\int
_{0}^1\tr\bigl[H(t)\bigr] \,\mathrm{d}t
\\
&&\quad\qquad{}+
\frac{n^{1-4\gamma}\mu_2^{2}(c_{\max}^5-c_{\min}^5)}{20}\int
_0^1\bigl[
\betav''(t)\bigr]^\top M(t)
\betav''(t) \,\mathrm{d}t \biggr\}\Rightarrow N
\bigl(0,\bigl(\sigma^*\bigr)^2\bigr),
\end{eqnarray*}
where
\begin{eqnarray*}
\bigl(\sigma^*\bigr)^2&=&\int_{\R}Q^2(c_{\max},t)
\,\mathrm{d}t\int_{0}^1\tr\bigl[H(t)^2
\bigr] \,\mathrm{d}t\quad\mbox{ and}\\
Q(x,y)&=&\int_{c_{\min}}^{x}
\bigl[2K(y/z)-K\ast K(y/z)\bigr]/z \,\mathrm{d}z.
\end{eqnarray*}
\end{theorem}

Now we consider the local power of $\lambda^*_n$ under the alternative
$H_{a,n}$ specified in Proposition~\ref{proplocalpower}. By Theorem
\ref{thmnewtest} and similar arguments as those of Proposition \ref
{proplocalpower}, it is easy to show that the asymptotic local power
of $\lambda^*_n$ with level $\alpha$
%
%e19 #&#
%
\begin{eqnarray}
\label{powernewtest} \beta^*_{\alpha}(c_{\min},c_{\max})=
\Phi(R_2-z_{1-\alpha})
\nonumber
\\[-8pt]
\\[-8pt]
\eqntext{\mbox{where } \displaystyle R_2=
\frac{(c_{\max}-c_{\min})F_1-(c^5_{\max}-c^5_{\min})\mu
_2^2F_2/20}{\sigma^*}.\qquad\qquad}
\end{eqnarray}
Suppose that $\lambda_n$ is asymptotically unbiased; namely $R_1>0$.
From \eqref{powernewtest} and \eqref{powerglrt}, we observe that
$\lambda^*_n$ is asymptotically more powerful than $\lambda_n$ if and
only if $R_2/R_1>1$. Simple calculations show that
\begin{eqnarray*}
R_2/R_1=\frac{[(c_{\max}-c_{\min})F_1-(c^5_{\max}-c^5_{\min})\mu
_2^2F_2/20]\sqrt{\int_{\R}\tilde{K}^2(t) \,\mathrm
{d}t}}{[c^{1/2}F_1-c^{9/2}\mu_2^2F_2/4]\sqrt{\int_{\R}Q^2(c_{\max
},t) \,\mathrm{d}t}}.
\end{eqnarray*}
An interesting observation from the above equation is that $R_2/R_1$
does not depend on the dependence or the non-stationarity structure of
the data. Furthermore, we have the following result.
%
%pr3 #&#
%
\begin{proposition}\label{proppower}
Under $H_{a,n}$ and the assumptions of Proposition \ref
{proplocalpower}, we have
%
%e20 #&#
%
\begin{equation}
\label{eqlargepower} \sup_{0<c_{\min}<c_{\max}<\infty}\beta
^*_{\alpha}(c_{\min},c_{\max})
\ge\sup_{0<c<\infty}\beta_\alpha(c).
\end{equation}
\end{proposition}
Proposition \ref{proppower} claims that, asymptotically, the averaged
test $\lambda^*_n$ is at least as powerful as the test which is based
on the maximum generalized likelihood ratio. The result is very general
in the sense that it does not depend on the nature of the local
alternative $\mathbf{f}_n(\cdot)$, the dependence structure of the
data or
the kernel function. When we restrict ourselves to a specific kernel
function, the power comparison can be more exact. Let us consider the
following example:

%f1 #&#
%
\begin{figure}[b]

\includegraphics{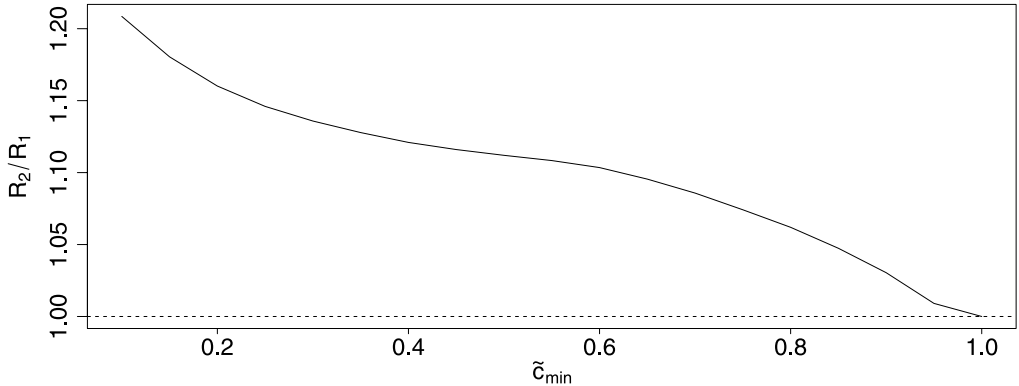}

\caption{Ratio $R_2/R_1$ as a function of $\tilde{c}_{\min}$
in Example \protect\ref{expower}. The uniform kernel is used.}\label
{fidpower}
\end{figure}

%ex5 #&#
%
\begin{example}\label{expower}
Suppose that $\lambda_n$ is asymptotically unbiased and that the
bandwidth for $\lambda_n$ is chosen as $cn^{-2/9}$. Let $c_{\min
}=\tilde
{c}_{\min}c$ for some fixed $\tilde{c}_{\min}\le1$ and let $c_{\max
}=\tilde{c}_{\max}c$ such that $\tilde{c}_{\max}$ solves the equation
$x^4+\tilde{c}_{\min}x^3+(\tilde{c}_{\min})^2x^2+(\tilde{c}_{\min
})^3x+(\tilde{c}_{\min})^4=5$. Choosing $c_{\max}$ in the latter way
insures that $F_1$ and $F_2$ do not enter the ratio $R_2/R_1$ and hence
the power comparison is relatively simple. Now simple calculations show that
%
%e21 #&#
%
\begin{eqnarray}
\label{eqpowercomparison} R_2/R_1=\frac{(\tilde{c}_{\max}-\tilde
{c}_{\min})\sqrt{\int_{\R}\tilde
{K}^2(t) \,\mathrm{d}t}}{\sqrt{\int_{\R} (\int_{\tilde{c}_{\min
}}^{\tilde
{c}_{\max}}[2K(y/z)-K\ast K(y/z)]/z \,\mathrm{d}z )^2 \,\mathrm{d}y}}.
\end{eqnarray}
An application of the Cauchy--Schwarz inequality similar to the proof of
Proposition \ref{proppower} shows that $\sup_{0<\tilde{c}_{\min}\le
1}R_2/R_1\ge1$ regardless of the kernel function. Now let us consider
the uniform kernel $K(x)=I\{|x|\le1\}/2$. Figure \ref{fidpower}\vadjust{\goodbreak}
shows $R_2/R_1$ as a function of $\tilde{c}_{\min}$. We observe from
the figure that the averaged test $\lambda^*_n$ is asymptotically more
powerful than $\lambda_n$ on $(0,1)$ regardless of the shape of the
alternative. %
Figure \ref{fidpower} further supports the use of the
averaged test.

%And the power performance of $\lambda^*_n$ is stable in this
%relatively wide interval. Indeed, the ratio $R_2/R_1$ ranges from 1 to
%1.02 in the latter range.
\end{example}

%Depending on the nature of the alternative and the kernel function,
%there are cases in which the equality in \eqref{eqlargepower} holds.

%s4.2 #&#
\subsection{Bias reduction and bandwidth range selection}\label{secbandwidth}
As we see from Theorem \ref{thmnewtest}, the asymptotic bias of
$\lambda^*_n$ involves the second derivative of $\betav(t)$ and the
estimation of the latter quantity is generally highly nontrivial.
Following the idea of Fan and Jiang \cite{FanJia07}, a prewhitening technique
can be used to alleviate the problem. More specifically, consider the
following null hypothesis:
\[
\tilde{H}_{0} \dvt\betav(\cdot)=\betav_{0}(\cdot,\theta)\qquad
\mbox{for some unknown } \theta_0\in\Omega\subset\R^q,
\]
where $\{\betav_{0}(\cdot,\theta)\dvt\theta\in\Omega\}$ is a parametric
family of smooth functions. Let $\hat{\theta}_0$ be a $\sqrt{n}$
consistent estimator of $\theta_0$ and define $\betav^*(t)=\betav
(\cdot
)-\betav_0(t,\hat{\theta}_0)$. Then by the similar arguments as those
of Proposition \ref{propequasemi}, the asymptotic bias and variance
of estimating $\theta_0$\vadjust{\goodbreak} is negligible in the current setting and hence
testing $\tilde{H}_0$ is equivalent to testing
\[
\breve{H}_{0} \dvt\betav^*(\cdot)=0\quad \mbox{versus}\quad
\breve{H}_{a} \dvt\betav^*(\cdot)\neq0.
\]
Then we can perform $\lambda^*_n$ to testing $\breve{H}_{0}$ with
transformed regression coefficients $\betav^*(\cdot)$ and response
$\breve{y}_i=y_i-\mathbf{x}_i^\top\betav_0(t,\hat{\theta}_0)$.
Note that
the local linear estimator of $\betav^*(\cdot)$ has no bias under
$\breve{H}_0$ and we can avoid the notorious problem of bias
estimation .

As mentioned in Fan and Jiang \cite{FanJia07}, a choice of larger bandwidth
favors smoother alternatives and a smaller bandwidth tends to detect
less smooth alternatives. Thanks to the introduction of the averaged
test, the sensitivity of the test to the choice of bandwidth is
alleviated due to the introduction of a group of bandwidths. On the
other hand, the correlation of $\lambda_n$ between nearby bandwidths
are usually quite high and hence in practice one only needs to average
the test over a grid of relatively separated bandwidths. Zhang \cite{Zha03}
found that the correlation between $\lambda_n(h)$ and $\lambda_n(ch)$
is quite high for $c=1.3$. As suggested by Fan and Jiang \cite
{FanJia07}, here
we recommend choosing the grid of three bandwidths $\tilde{b}_n/1.5$,
$\tilde{b}_n$ and $\tilde{b}_n\times1.5$ to represent small, medium
and large bandwidths and average the test over the latter grid. Here
$\tilde{b}_n=b_n^*\times n^{-1/45}$ and $b_n^*$ is the optimal
bandwidth for nonparametric curve estimation.

%s4.3 #&#
\subsection{The robust wild bootstrap}\label{secwildboots}
A direct implementation of the asymptotic distribution in Theorem \ref
{thmnewtest} may not perform satisfactorily in practice due to the
following two reasons. First, the convergence rate of test statistic
equals $\mathrm{O}(n^{-1/9})$ when bandwidth $b_n$ is chosen optimally. The rate
is quite slow and hence the asymptotic approximation may not be
accurate for moderate samples. Second, as we can see from the proof of
Lemma \ref{lemasy5} in Section \ref{secproofs}, the asymptotic normal
approximation is particularly rough at the boundaries of the time
interval for finite samples. As a remedy, we observe the following proposition.

%pr4 #&#
%
\begin{proposition}\label{propwildboots}
Let the bandwidth range be $[c_{\min}n^{-\gamma}, c_{\max}n^{-\gamma}]$
for some $0<c_{\min}<c_{\max}<\infty$. Suppose that either (1):
$\betav_0(\cdot)$ is a linear function or (2): $\gamma>2/9$. Then
under $H_0$,
condition~\textup{(A)} and the assumption that $\gamma<1/4$, on a possibly
richer probability space, there exist i.i.d. $p$-dimensional standard
Gaussian random vectors $V_1,\ldots,V_n$, such that
%
%e22 #&#
%
\begin{equation}
\label{eqwildboots} \lambda_n^*=\Phi_n+\mathrm{o}_\p
\bigl(\sqrt{n^{\gamma}}\bigr),
\end{equation}
where
\begin{eqnarray*}
\Phi_n=\int_{c_{\min}}^{c_{\max}} \Biggl\{2\sum
_{i=1}^n\tilde{V}^{\top
}_i
\bigl[\E\mathbf{ S}_{n,n(s)}(t_i)\bigr]^{-1}\tilde{
\mathbf{T}}_{n,n(s)}(t_i)-\sum_{i=1}^n
\bigl[\mathbf{ z}^{\top}_i\bigl[\E\mathbf{
S}_{n,n(s)}(t_i)\bigr]^{-1}\tilde{
\mathbf{T}}_{n,n(s)}(t_i)\bigr]^2 \Biggr\} \,
\mathrm{d}s
\end{eqnarray*}
with $n(s)=sn^{-\gamma}$, $\mathbf{z}_i=(\mathbf{x}_i^{\top
},\mathbf{0}_p^{\top
})^{\top}$, $\tilde{V}_i=(V^{\top}_i\Lambda^{1/2}(t_i),\mathbf
{0}_p^{\top
})^{\top}$, $\tilde{\mathbf{T}}_{n,b}(t)=(\tilde{\mathbf{T}}^{\top
}_{n,0,b}(t),\tilde{\mathbf{T}}^{\top}_{n,1,b}(t))^{\top}$ and
%
%e23 #&#
%
\begin{eqnarray}
\label{eqTn} \tilde{\mathbf{T}}_{n,l,b}(t)=(nb)^{-1}\sum
_{i=1}^n\Lambda^{1/2}(t_i)V_i
\bigl[(t_i-t)/b\bigr]^l K_{b}(t_i-t),\qquad
l=0,1.
\end{eqnarray}
\end{proposition}
Proposition \ref{propwildboots} follows easily from \eqref{eq1} and
Lemma \ref{lemasy3} in Section \ref{secproofs}. Details are omitted.
The latter proposition claims that $\lambda^*_n$ can be well
approximated by a Gaussian quadratic form $\Phi_n$. In particular, we
observe from the proofs in Section \ref{secproofs} that the
approximation is accurate at the boundaries due to the fact that it
directly mimics the form of the test statistic. When implementing~$\lambda_n^*$,
we recommend generating a large (say of size 1000)
sample of i.i.d. copies of $\Phi_n$ and use the resulting empirical
distribution to approximate that of $\lambda_n^*$ under the null
hypothesis and obtain the $p$-value of the test.

As we suggested in Section \ref{secbandwidth}, in practice, one
usually uses a grid of bandwidths $\mathcal{ B}=\{c_{\min}n^{-\gamma
}=b_1<b_2<\cdots<b_M=c_{\max}n^{-\gamma}\}$ and calculate $\lambda
_n^*(\mathcal{ B})=\sum_{i=1}^M(\RSS_{0}-\RSS_{a}(b_i))$. To
perform wild
bootstrap in those cases, one compares $\lambda_n^*(\mathcal{ B})$ to the
simulated quantiles of
\begin{eqnarray*}
\Phi_n(\mathcal{ B}):=\sum_{j=1}^M
\Biggl\{2\sum_{i=1}^n\tilde{V}^{\top}_i
\bigl[\E\mathbf{ S}_{n,b_j}(t_i)\bigr]^{-1}\tilde{
\mathbf{T}}_{n,b_j}(t_i)-\sum_{i=1}^n
\bigl[\mathbf{ z}^{\top}_i\bigl[\E\mathbf{
S}_{n,b_j}(t_i)\bigr]^{-1}\tilde{
\mathbf{T}}_{n,b_j}(t_i)\bigr]^2 \Biggr\}
\end{eqnarray*}
to calculate the $p$-value of the test. In Section \ref{secsimu}, we
shall conduct a simulation study to compare the finite sample
performance of the wild bootstrap and the direct implementation of the
asymptotic distribution.
%w(b_jn^{\alpha})

If one is interested in the semiparametric testing problem $H_{01}$
versus $H_{a1}$ in \eqref{eqseminull}, then the corresponding
averaged test is

%e24 #&#
%
\begin{equation}
\label{eqtestcomplex} \lambda^*_{1n}=\int_{c_{\min}}^{c_{\max}}
\bigl(\RSS_{1}\bigl(zn^{-\gamma}\bigr)-\RSS_{a}
\bigl(zn^{-\gamma}\bigr)\bigr) \,\mathrm{d}z.
\end{equation}

Write $\varepsilon_i=([\varepsilon^{(1)}_{i}]^{\top},[\varepsilon
^{(2)}_{i}]^\top)^{\top}$ and $V_i=([V^{(1)}_{i}]^{\top
},[V^{(2)}_{i}]^\top)^{\top}$, where $\varepsilon^{(1)}_{i}$ and
$V^{(1)}_{i}$ are $p_1$ dimensional. Define $\mathbf{S}^{(2)}_{n,b}$,
$\mathbf{S}^{(2)}_{n,l,b}$, $\mathbf{z}_{i}^{(2)}$, $\tilde{V}^{(2)}_i$,
$\mathbf{T}^{(2)}_{n} \mathbf{T}^{(2)}_{nl}$, $\tilde{\mathbf
{T}}^{(2)}_{n}$,
$\tilde{\mathbf{T}}^{(2)}_{nl}$ and $\Phi^{(2)}_n$ in the same way as
their counterparts without the superscript $^{(2)}$ with $\mathbf{x}_i$,
$\varepsilon_i$, $\Lambda(t)$ and $V_i$ therein replaced by $\mathbf{
x}^{(2)}_i$, $\varepsilon^{(2)}_i$, $\Lambda_{22}(t)$ and $V^{(2)}_i$,
respectively. We have the following
proposition.
%
%pr5 #&#
%
\begin{proposition}\label{propwildbootscomplex}
Suppose that $1/4>\gamma>2/9$. Then under $H_{01}$ and condition \textup{(A)},
on a possibly richer probability space, there exist i.i.d.
$p$-dimensional standard Gaussian random vectors $V_1,\ldots,V_n$,
such that
%
%e25 #&#
%
\begin{equation}
\label{eqwildbootscomplex} \lambda_{1n}^*=\Phi_n-
\Phi^{(2)}_n+\mathrm{o}_\p\bigl(\sqrt{n^{\gamma}}
\bigr).
\end{equation}
\end{proposition}
Note that $\Phi_n-\Phi^{(2)}_n$ is a quadratic form of $V_1,\ldots
,V_n$. By Proposition \ref{propwildbootscomplex}, in practice, one
could generate a large sample of i.i.d. copies of $\Phi_n-\Phi^{(2)}_n$
to obtain the $p$-value of testing $H_{01}$.

%s4.4 #&#
\subsection{Long-run covariance matrix estimation}
By Lemma \ref{lemsninv} in Section \ref{secproofs}, $\E\mathbf{
S}_{n,n(s)}(t_i)$ in Proposition \ref{propwildboots} can be well
approximated by $\mathbf{
S}_{n,n(s)}(t_i)$. Therefore, in order to implement the wild bootstrap,
one only needs to estimate the long-run covariance matrix $\Lambda
(\cdot
)$. Here we suggest using the local lag window estimate of $\Lambda
(\cdot)$ proposed in Zhou and Wu \cite{ZhoWu10}. For the sake of completeness,
we will briefly introduce the estimator here. We refer to the latter
paper for more details including the derivation of convergence rates of
the estimator and the choice of smoothing parameters.

Define $\hat{\mathbf{L}}_i:=\mathbf{x}_i\hat{\varepsilon}_i$,
where $\hat
{\varepsilon}_i$'s are the residuals under the alternative. For a
window size $m$ and a bandwidth $\tau_n$, $\Lambda(\cdot)$ can be
estimated by
\[
\hat{\Lambda}(\cdot)=\sum_{i=1}^n\omega(
\cdot,i)\Delta_i \qquad\mbox{where } \omega(\cdot,i)=\frac{K_{\tau
_n}(t_i-\cdot)}{\sum_{j=1}^nK_{\tau
_n}(t_j-\cdot)}
\]
and $\Delta_i=(\sum_{j=-m}^m\hat{\mathbf{L}}_{i+j})(\sum
_{j=-m}^m\hat{\mathbf{
L}}^{\top}_{i+j})/(2m+1)$. Zhou and Wu \cite{ZhoWu10} showed that
$\hat{\Lambda
}(t)$ is always positive semidefinite and has convergence rate
$\mathrm{O}(n^{-2/7})$ when $m=\mathrm{O}(n^{2/7})$ and $\tau_n=\mathrm{O}(n^{-1/7})$.

%s5 #&#
\section{Simulation studies}\label{secsimu}
In this section, we shall design simulations to study the accuracy of
the wild bootstrap procedure of the paper and compare it with that of
the bootstrap procedure of Fan and Jiang \cite{FanJia07} and the
method of
direct implementation of the asymptotic distribution in \eqref
{eqnull}. Let us consider the following model
%
%e26 #&#
%
\begin{equation}
\label{modelsimu} y_i=\beta_1(t_i)+
\beta_2(t_i)x_{2i}+\varepsilon_i
\end{equation}
and the test $H_0\dvt \beta_1(\cdot)=\beta_2(\cdot)=0$. The following
four scenarios are considered in order to investigate the effects of
endogeneity, non-stationarity and temporal dependence.
\begin{longlist}[Scenario (a)]
\item[Scenario (a).] In this case $x_{2i}$'s are i.i.d.
exponential random variables with mean 1 and $\varepsilon_i$'s are
i.i.d. standard normal. The two processes $\{x_{2i}\}$ and $\{
\varepsilon_i\}$ are independent. The latter design satisfies the
conditions in Fan, Zhang and Zhang \cite{FanZhaZha01} and hence it is
expected that
the bootstrap procedure in Fan and Jiang \cite{FanJia07} will work in
this case.

\item[Scenario (b).] In this scenario $x_{2i}$'s are i.i.d.
exponential random variables with mean 1 and $\varepsilon
_i=x_{2i}\zeta_i$, where $\zeta_i$'s are i.i.d. standard normal and
are independent
of $\{x_{2i}\}$. In scenario (b) we are interested in investigating the
effect of endogeneity on the behavior of GLRT.

\item[Scenario (c).] Let $x_{2i}$'s be independent student $t$
random variables and the degrees of freedom of $x_{2i}=5+10t_i$. Let
$\varepsilon_i=\exp(-1/t_i)/(100t_i^4)\zeta_i$, where $\zeta_i$'s are
i.i.d. standard normal. Further let $x_{2i}$'s and $\varepsilon_i$'s be
independent. Note that $\{\varepsilon_i\}$ is a locally stationary
process with time-varying variance and $\{x_{2i}\}$ is locally
stationary process with smoothly varying tail index. In this case, we
are investigating the effect of non-stationarity on the behavior of GLRT.

\item[Scenario (d).] Let $x_{2i}=\epsilon_i\epsilon_{i-1}$,
where $\epsilon_i$'s are i.i.d. standard normal. Let $\varepsilon
_i=0.5\varepsilon_{i-1}+\zeta_{i}$, where $\zeta_i$'s are i.i.d.
standard normal. Further let $\{\epsilon_i\}$ be independent of $\{
\zeta_i\}$. Note $\{x_{2i}\}$ and $\{\varepsilon_i\}$ are two stationary
weakly dependent processes. In this case we are interested in
investigating the effect of temporal dependence on the behavior of GLRT.
\end{longlist}
We consider two different sample sizes, $n=200$ and $400$. We compare
three different methods, namely the robust wild bootstrap test \eqref
{eqwildboots} (WILD), test based on the asymptotic distribution \eqref
{eqnull} (ASYM) and the residual bootstrap test of Fan and Jiang
\cite{FanJia07} (IID). Both the single bandwidth test $\lambda_n$ in
\eqref
{eqglrtoriginal} and the suggested averaged test $\lambda_n^*$ in
\eqref{eqtest} are considered. For the averaged test, the bandwidth
ranges are selected as $[\tilde{b}_n/1.5,1.5\tilde{b}_n]$ according to
the discussion in Section \ref{secbandwidth}. To investigate the
sensitivity of the accuracy of the wild bootstrap method on the choice
of bandwidth, three different bandwidths, namely $0.15, 0.25$ and
$0.35$ are considered in the simulation. Based on 500 replications, the
simulated type I error rates at $10\%$ nominal level are summarized in
Table \ref{tab1} below.\looseness=1

\begin{table}
\def\arraystretch{0.9}
\caption{Simulated type I error rates (in percentage) for the
wild bootstrap test \protect\eqref{eqwildboots} (WILD), test based on the
asymptotic distribution \protect\eqref{eqnull} (ASYM) and the
bootstrap test
of Fan and Jiang \cite{FanJia07} (IID) with nominal level $10\%$ under
scenarios
\textup{(a)}, \textup{(b)}, \textup{(c)} and \textup{(d)}. For the
averaged test $\lambda_n^*$, the
bandwidth range is selected as $[\tilde{b}_n/1.5, 1.5\tilde{b}_n]$.
Series length $n=200$ and $400$ with $500$ replicates}\label{tab1}
\begin{tabular*}{\textwidth}{@{\extracolsep{\fill
}}lld{2.1}d{2.2}d{2.1}d{2.1}d{2.1}d{2.1}d{2.1}d{2.1}@{}}
\hline
& & \multicolumn{4}{l}{$n=200$} &\multicolumn{4}{l@{}}{$n=400$}\\[-6pt]
& & \multicolumn{4}{l}{\hrulefill} &\multicolumn{4}{l@{}}{\hrulefill
}\\
\multicolumn{1}{@{}l}{Method} & & \multicolumn{1}{l}{(a)} &
\multicolumn{1}{l}{(b)}& \multicolumn{1}{l}{(c)} & \multicolumn
{1}{l}{(d)} &
\multicolumn{1}{l}{(a)} & \multicolumn{1}{l}{(b)} & \multicolumn
{1}{l}{(c)} &
\multicolumn{1}{l@{}}{(d)} \\
\hline
\multicolumn{10}{@{}l}{Averaged test $\lambda_n^*$}\\
WILD & $\tilde{b}_n=0.15$&
7.5 & 7.4 & 10.4 & 7.1 & 8.1 & 8
& 9.7 & 9.1 \\
WILD & $\tilde{b}_n=0.25$&
8.5 & 8.15 & 10.2 & 7.7 & 8.5 &
8.4 & 9.8 & 9.7\\
WILD & $\tilde{b}_n=0.35$&
8.9 & 8.7 & 10 & 7.7 & 8.7 & 9.1
& 9.2 & 9.8\\
ASYM & $\tilde{b}_n=0.15$&
35.4 & 14.4 & 18.8 & 28.2 & 38.3 &
18.5 & 15.0 & 33\\
ASYM & $\tilde{b}_n=0.25$&
39.1 & 18.5 & 19.4 & 33.3 & 39.9 &
21.2 & 17.8 & 36.3 \\
ASYM & $\tilde{b}_n=0.35$&
44.1 & 21.4 & 18.0 & 36.2 & 44.5 &
23.8 & 20.7 & 38.4 \\
IID & $\tilde{b}_n=0.15$&
10.4 & 83.6 & 20.5 & 68.8 & 11.9 &
87.7 & 15.7 & 73.3\\
IID & $\tilde{b}_n=0.25$&
11.4 & 79.6 & 19.1 & 61.9 & 9.9 &
82.7 & 17.9 & 63.8 \\
IID & $\tilde{b}_n=0.35$&
11.0 & 74.3 & 17.8 & 55.9 & 10.2 &
78.8 & 19.8 & 56.8\\[3pt]
\multicolumn{10}{@{}l}{Single bandwidth test $\lambda_n$}\\
WILD & $b_n=0.15$&
5.0 & 5.8 & 10.2 & 5.8 & 8.6 & 7.2 &
11.2 & 9.4 \\
WILD & $b_n=0.25$&
8.2 & 7.8 & 9.4 & 8.8 & 9.2 & 8.2 &
10.2 & 11.6\\
WILD & $b_n=0.35$&
9.8 & 9.2 & 9.0 & 8.2 & 11.2 & 9.6 &
11.2 & 11.4\\
ASYM & $b_n=0.15$&
32.2 & 13.2 & 17.8 & 27.8 & 27.4 & 16.8 &
13.8 & 30\\
ASYM & $b_n=0.25$&
36.2 & 19.6 & 20.4 & 36.8 & 29 & 20.2 &
16.8 & 36.6 \\
ASYM & $b_n=0.35$&
43.6 & 21.2 & 20.4 & 38.8 & 34 & 22 &
18 & 38 \\
IID & $b_n=0.15$&
8.2 & 86.8 & 20.8 & 73.2 & 10.8 & 89 &
15.2 & 76.2\\
IID & $b_n=0.25$&
7.8 & 82.2 & 20.6 & 63 & 9.4 & 80.2 & 18
& 63.4 \\
IID & $b_n=0.35$&
10.4 & 76.2 & 17.4 & 55.6 & 12 & 77.2 &
17.8 & 56.6\\
\hline
\end{tabular*}\vspace*{-3pt}
\end{table}

We observe from Table \ref{tab1} that, for the robust wild bootstrap, the
simulated type I errors of the averaged test and the single bandwidth
test are reasonably close to the nominal and the performance is stable
for all four cases when $n=400$. For $n=200$, the robust bootstrap is
slightly anti-conservative in cases (a), (b) and (d) for small
bandwidths. As we expected, the averaged test performs more stably than
the single bandwidth test. On the other hand, we observe that tests
based on the asymptotic distribution do not perform well for moderately
large samples. As we discussed in Section \ref{secwildboots}, the
reason is due to the slow convergence of the test statistic and the
rough approximation of the asymptotic distribution at the boundaries.
The residual wild bootstrap performs slightly better than our robust
wild bootstrap for i.i.d. data without endogeneity. However, we observe
that the residual bootstrap is no longer consistent under
non-stationarity, temporal dependence or endogeneity, which is
consistent with our theoretical findings.\looseness=1

%s6 #&#
\section{Proofs}\label{secproofs}

Note that under the null hypothesis $H_0$,
%
%e27 #&#
%
\begin{eqnarray}
\label{eqp1} \RSS_a-\RSS_0=2\sum
_{i=1}^n\mathbf{ x}^{\top}_i
\varepsilon_i\bigl(\betav(t_i)-\hat{
\betav}(t_i)\bigr)+\sum_{i=1}^n
\bigl\{ \mathbf{ x}^{\top}_i\bigl(\betav(t_i)-
\hat{\betav}(t_i)\bigr)\bigr\}^2:=2I_n+\mathit{II}_n.
\end{eqnarray}
On the other hand, by \eqref{eqsol},
%
%e28 #&#
%
\begin{equation}
\label{eq6} \label{eqthm11} \mathbf{S}_n(t) \bigl(\hat{\etav}(t)-
\etav(t)\bigr) =\pmatrix{b_n^2
\mathbf{ S}_{n,2}(t) \bigl(\betav''(t)+\mathrm{o}(1)
\bigr)/2
\vspace*{2pt}\cr
b_n^2\mathbf{ S}_{n,3}(t) \bigl(
\betav''(t)+\mathrm{o}(1)\bigr)/2
}
 +\pmatrix{\mathbf{
T}_{n,0}(t)
\vspace*{2pt}\cr
\mathbf{T}_{n,1}(t)
}
:=\mathbf{B}_n(t)+\mathbf{T}_n(t),
\end{equation}
where $\etav(t)=(\betav^{\top}(t),b_n\betav'^{\top}(t))^{\top}$, and
\begin{eqnarray*}
\mathbf{T}_{n,l}(t)=r_n^2\sum
_{i=1}^n\mathbf{x}_i
\varepsilon_i \bigl[(t_i-t)/b_n
\bigr]^l K_{b_n}(t_i-t),\qquad  l=0,1, \ldots
\end{eqnarray*}
with $r_n:=1/\sqrt{nb_n}$. In \eqref{eqthm11}, $\mathbf{B}_n(t)$
corresponds to the bias of the
local linear estimate at time $t$. Lemmas \ref{lembias1} and
\ref{lembias2} below control the asymptotic influence of the bias
term $\mathbf{B}_n(\cdot)$ on $\RSS_a-\RSS_0$.

%le1 #&#
%
\begin{lemma}\label{lembias1}
Define $\mathbf{z}_i=(\mathbf{x}_i^{\top},\mathbf{0}_p^{\top
})^{\top}$, where
$\mathbf{0}_p$ is the column vector of $p$ zeros. Under condition~\textup{(A)},
we have $-I_n=D_{n1}+\mathrm{O}_\p(\sqrt{n}b_n^2)$,
where $D_{n1} :=\sum_{i=1}^n\mathbf{z}^{\top}_i\varepsilon_i\mathbf{
S}^{-1}_n(t_i)\mathbf{T}_n(t_i)$.
\end{lemma}
\begin{pf}
By \eqref{eqp1} and \eqref{eqthm11}, we have
\[
-I_n-D_{n1}=\sum_{i=1}^n
\mathbf{z}^{\top}_i\varepsilon_i\mathbf{
S}^{-1}_n(t_i)\mathbf{B}_n(t_i).
\]
Define $ID_{n1}=\E[(-I_n-D_{n1})^2|\FF_n]$ and $\PP_i(\cdot)=\E
(\cdot
|\GG_i)-\E(\cdot|\GG_{i-1}).$
Recall that $\GG_i=(\ldots,\xi_{i-1},\xi_i)$. Using the facts that
$H(t, \GG_i)=\sum_{j=-\infty}^i\PP_jH(t, \GG_i)$ and $\PP_i$ and
$\PP_j$ are orthogonal for $i\neq j$, elementary calculations show
that
\begin{eqnarray*}
ID_{n1}&=&\sum_{i=1}^n\sum
_{j=1}^n\sum_{k=-\infty}^n
\E\bigl[\PP_kH(t_i,\GG_i)
\PP_kH(t_j,\GG_j)\bigr]\\[-2pt]
&&\hspace*{54pt}{}\times V(t_i,
\FF_i)\mathbf{ S}^{-1}_n(t_i)
\mathbf{B}_n(t_i)V(t_j,\FF_j)
\mathbf{S}^{-1}_n(t_j)\mathbf{
B}_n(t_j).
\end{eqnarray*}
Let $\delta_H(k,p)=0$ if $k<0$. Note that
\begin{eqnarray*}
\sum_{k=-\infty}^n\bigl|\E\bigl[
\PP_kH(t_i,\GG_i)\PP_kH(t_j,
\GG_j)\bigr]\bigr|&\le& \sum_{k=-\infty}^n
\bigl\|\PP_kH(t_i,\GG_i)\bigr\|\bigl\|\PP_kH(t_j,
\GG_j)\bigr\|
\\[-2pt]
&\le&\sum_{k=-\infty}^n\delta_H(i-k,2)
\delta_H(j-k,2)
\\[-2pt]
%&\le&C\sum_{k=-\infty}^{\min(i,j)}(i-k+1)^{-2}(j-k+1)^{-2}\\
&\le& C\bigl(|i-j|+1\bigr)^{-2}.
\end{eqnarray*}
On the other hand, by Lemma \ref{lemsninv}, the H\"{o}lder's
inequality and similar arguments as those of Lemma 6 in Zhou and Wu
\cite{ZhoWu10}, we have
\begin{eqnarray*}
&&\E\bigl|V(t_i,\FF_i)\mathbf{S}^{-1}_n(t_i)
\mathbf{B}_n(t_i)V(t_j,\FF_j)
\mathbf{ S}^{-1}_n(t_j)\mathbf{B}_n(t_j)\bigr|
\\[-2pt]
&&\quad\le\bigl\|V(t_i,\FF_i)\bigr\|_{4}\bigl\|
\mathbf{S}^{-1}_n(t_i)\bigr\|_{8}\bigl\|
\mathbf{ B}_n(t_i)\bigr\|_8\bigl\|V(t_j,
\FF_j)\bigr\|_{4}\bigl\|\mathbf{S}^{-1}_n(t_j)
\bigr\|_8\bigl\|\mathbf{ B}_n(t_j)\bigr\|_8\le
Cb_n^4.
\end{eqnarray*}
Therefore, $\E ID_{n1}\le C\sum_{i=1}^n\sum
_{j=1}^n(|i-j|+1)^{-2}b_n^4\le
Cnb_n^4.$
Note that $\E(-I_n-D_{n1})^2=\E ID_{n1}$. Therefore, this lemma
follows.
\end{pf}

%le2 #&#
%
\begin{lemma}\label{lembias2}
Under condition \textup{(A)} and the assumption that $nb_n^{5/2}\rightarrow
\infty
$, we have
\[
\mathit{II}_n=D_{n2}+\frac{nb_n^4\mu_2^{2}}{4}\int_0^1
\bigl[\betav''(t)\bigr]^\top M(t)
\betav''(t) \,\mathrm{d}t+\mathrm{o}_\p
\bigl(nb_n^4\bigr),
\]
where $D_{n2} :=\sum_{i=1}^n\{\mathbf{z}^{\top}_i\mathbf{
S}^{-1}_n(t_i)\mathbf{T}_n(t_i)\}^2$.\vadjust{\goodbreak}
\end{lemma}
\begin{pf}
By \eqref{eqp1} and \eqref{eqthm11}, we have
\begin{eqnarray*}
\label{eqp2} \mathit{II}_n-D_{n2}&=&\sum
_{i=1}^n\bigl(\mathbf{z}^{\top}_i
\mathbf{S}^{-1}_n(t_i)\mathbf{
B}_n(t_i)\bigr)^2+2\sum
_{i=1}^n\mathbf{z}^{\top}_i
\mathbf{S}^{-1}_n(t_i)\mathbf{
B}_n(t_i)\mathbf{z}^{\top}_i
\mathbf{S}^{-1}_n(t_i)\mathbf{T}_n(t_i)
\\
&:=&I D^*_{n2}+2I D^{**}_{n2}.
\end{eqnarray*}
By Lemma \ref{lemsninv} and the H\"{o}lder's inequality, it follows that
\[
I D^*_{n2}-\sum_{i=1}^n\bigl\{
\mathbf{z}^{\top}_i\bigl[\E\mathbf{ S}_n(t_i)
\bigr]^{-1}\mathbf{B}_n(t_i)\bigr
\}^2=\mathrm{O}_\p\bigl(nb_n^4/
\sqrt{nb_n}\bigr).
\]
By condition (A4) and the similar arguments as those in the proof of
Lemma \ref{lembias1}, we have
\begin{eqnarray*}
\sum_{i=1}^n\bigl\{\mathbf{z}^{\top}_i
\bigl[\E\mathbf{ S}_n(t_i)\bigr]^{-1}
\mathbf{B}_n(t_i)\bigr\}^2-\E\Biggl[\sum
_{i=1}^n\bigl\{\mathbf{z}^{\top}_i
\bigl[\E\mathbf{ S}_n(t_i)\bigr]^{-1}
\mathbf{B}_n(t_i)\bigr\}^2
\Biggr]=\mathrm{O}_\p\bigl(\sqrt{n}b_n^4\bigr).
\end{eqnarray*}
It is easy to see that, for $i=1,2,\ldots, n$,
\[
\E\bigl(\mathbf{z}^{\top}_i\bigl[\E\mathbf{S}_n(t_i)
\bigr]^{-1}\mathbf{B}_n(t_i)
\bigr)^2-b_n^4\E\biggl(\mathbf{z}^{\top}_i
\bigl[\E\mathbf{S}_n(t_i)\bigr]^{-1}\pmatrix{\mathbf{ S}_{n,2}(t_i)
\betav''(t_i)/2
\vspace*{2pt}\cr
\mathbf{ S}_{n,3}(t_i)\betav''(t_i)/2
}
 \biggr)^2=\mathrm{o}\bigl(b_n^4
\bigr).
\]
Additionally, by Lemma \ref{lemsninv} and simple algebra, we have
\begin{eqnarray*}
\sum_{i=1}^n\E\biggl(
\mathbf{z}^{\top}_i\bigl[\E\mathbf{S}_n(t_i)
\bigr]^{-1}\pmatrix{\mathbf{
S}_{n,2}(t_i)\betav''(t_i)/2
\vspace*{2pt}\cr
\mathbf{ S}_{n,3}(t_i)\betav''(t_i)/2
}
 \biggr)^2=n\mu^2_2
\int_0^1\bigl[\betav''(t)
\bigr]^\top M(t)\betav''(t) \,
\mathrm{d}t/4+\mathrm{o}(n).
\end{eqnarray*}
Therefore, $I D^*_{n2}=nb_n^4\mu^2_2\int_0^1[\betav''(t)]^\top
M(t)\betav''(t) \,\mathrm{d}t/4+\mathrm{o}_p(nb_n^4)$. Furthermore,
\begin{eqnarray*}
I D^{**}_{n2}=r_n^2\sum
_{j=1}^n\sum_{i=1}^n
\mathbf{z}^{\top}_i\mathbf{ S}^{-1}_n(t_i)
\mathbf{B}_n(t_i)\mathbf{z}^{\top}_i
\mathbf{S}^{-1}_n(t_i)\mathbf{
x}_jK_{b_n}(t_i-t_j)
\varepsilon_j.
\end{eqnarray*}
Recall that $r_n=1/\sqrt{nb_n}$. Following the similar arguments as those
in the proof of Lemma \ref{lembias1}, we have
$I D^{**}_{n2}=\mathrm{O}_\p(\sqrt{b^3_n})=\mathrm{o}_\p(nb_n^4)$. Details are omitted.
Hence, the lemma follows.
\end{pf}

%Lemmas \ref{lembias1} and \ref{lembias2} conclude that, for GLRT
%with the local linear kernel estimate $\hat{\betav}(\cdot)$, the
%asymptotic influence of the bias term $\mathbf{B}(\cdot)$ is of the
%order $\mathrm{O}_\p(nb_n^4)$, which is negligible under the bandwidth
%constraint $b_n=o(n^{-2/9})$.

%le3 #&#
%
\begin{lemma}\label{lemasy1}
Under condition \textup{(A)} and the assumption that $nb_n^{3}\rightarrow\infty
$, we have
\[
D_{n1}=\bar{D}_{n1}+\mathrm{o}_\p(1/
\sqrt{b_n}),
\]
where $\bar{D}_{n1}=\sum_{i=1}^n\mathbf{z}^{\top}_i\varepsilon
_i[\E\mathbf{
S}_n(t_i)]^{-1}\mathbf{T}_n(t_i)$.
\end{lemma}
\begin{pf}
Let $\barid_{n1}=D_{n1}-\bar{D}_{n1}$ and $\mathbf{IS}_n(t)=\mathbf{
S}_n^{-1}(t)-[\E\mathbf{S}_n(t)]^{-1}$. Then
\[
\barid_{n1}=\sum_{i=1}^n
\mathbf{z}^{\top}_i\varepsilon_i\mathbf{
IS}_n(t_i)\mathbf{T}_n(t_i).
\]
Let $\mathbf{A}_{n,k}=\sum_{i=1}^k\mathbf{z}^{\top}_i\varepsilon
_i\mathbf{
IS}_n(t_i)$ and $\mathbf{A}_{n,0}=0$. Then by Lemma \ref{lemsninv} and
the similar arguments as those of Lemma \ref{lembias1}, it is easy
to show that $\max_{1\le k\le n}\|\mathbf{A}_{n,k}\|_4\le Cr_n\sqrt{n}$.
Note that
\[
\barid_{n1}=\sum_{i=1}^n(
\mathbf{A}_{n,i}-\mathbf{A}_{n,i-1})\mathbf{
T}_n(t_i)=\sum_{i=1}^{n-1}
\mathbf{A}_{n,i}\bigl(\mathbf{T}_n(t_i)-
\mathbf{ T}_n(t_{i-1})\bigr)+\mathbf{A}_{n,n}
\mathbf{T}_n(t_n).
\]
By the similar arguments as those of Lemma \ref{lembias1}, we have
%
%e29 #&#
%
\begin{equation}
\label{eq4} \max_{1\le i\le n}\bigl\|\mathbf{T}_n(t_i)-
\mathbf{ T}_n(t_{i-1})\bigr\|_4\le
Cr_n^3
\end{equation}
and $\|\mathbf{T}_n(t_n)\|_4=\mathrm{O}(r_n)$.
Therefore,
\begin{eqnarray*}
\|\barid_{n1}\|&\le& \sum_{i=1}^{n-1}
\|\mathbf{A}_{n,i}\|_4\bigl\|\mathbf{ T}_n(t_i)-
\mathbf{T}_n(t_{i-1})\bigr\|_4+\|\mathbf{A}_{n,n}
\|_4\bigl\|\mathbf{ T}_n(t_n)\bigr\|_4
\\
&\le&C\Biggl(\sum_{i=1}^{n-1}r_n
\sqrt{n}r_n^3+r_n\sqrt{n}r_n
\Biggr)=\mathrm{O}\bigl(1/\bigl(\sqrt{n}b_n^2\bigr)\bigr)=\mathrm{o}(1/
\sqrt{b_n}).
\end{eqnarray*}
Therefore, the lemma follows.
\end{pf}

%IS}_n(t_i)\mathbf{z}^*_{ij}K_{b_n}(t_j-t_i)\varepsilon_i\varepsilon_j\\
%&=&r_n^2\sum_{i=1}^n\sum_{j=1}^n\mathbf{z}^{\top}_i\mathbf{IS}_n(t_i)
%z}^*_{ij}K_{b_n}(t_j-t_i)V(t_i, \FF_i)V(t_j, \FF_j)H(t_i,
%&:=&r_n^2\sum_{i=1}^n\sum_{j=1}^nA_n(i,j)H(t_i, \GG_i)H(t_j, \GG_j),
%where $\mathbf{z}^*_{ij}=(\mathbf{x}^{\top}_j,\mathbf{
%x}^{\top}_j(t_j-t_i)/b_n)^{\top}$. Recall that $r_n=1/\sqrt{nb_n}$.
%Note that $A_n(i,j)$ is measurable with respect to $\FF_n$ and
%$\GG_i$ is independent of $\FF_n$ for all $i$.

%Let $\barid_{n1}(\FF_n)=\E(\barid_{n1}|\FF_n)$ and
%$\barid^2_{n1}(\FF_n)=\var(\barid_{n1}|\FF_n)$. Then
%Since $\PP_i$ and $\PP_j$ are orthogonal if $i\neq j$,
%|\E[H(t_i, \GG_i)H(t_j,
%&=&|\sum_{k=-\infty}^n\E[\PP_kH(t_i,\GG_i)\PP_kH(t_j,\GG_j)]|\le
%&\le&\sum_{k=-\infty}^n\delta_H(i-k,2)\delta_H(j-k,2)\le
%C(|i-j|+1)^{-2}.
%On the other hand, by Lemma \ref{lemsninv} and H\"{o}lder's
%inequality, it is easy to show that
%Plugging the above results into \eqref{eqp3}, we have
%|\E\barid_{n1}|=|\E\barid_{n1}(\FF_n)|\le
%r_n^3\sum_{i}\sum_j(|i-j|+1)^{-2}K_{b_n}(t_j-t_i)\le Cr_n^3n.
%According to its definition, $\barid^2_{n1}(\FF_n)$ equals
%&&r_n^4\sum_{i=1}^n\sum_{j=1}^n\sum_{r=1}^n\sum_{m=1}^nA_n(i,j)A_n(r,m)
%&&=r_n^4\sum_{i=1}^n\sum_{j=1}^n\sum_{r=1}^n

%le4 #&#
%
\begin{lemma}\label{lemasy2}
Under condition \textup{(A)} and the assumption that $nb_n^3\rightarrow\infty$,
we have
\[
D_{n2}=\bar{D}_{n2}+\mathrm{o}_\p(1/
\sqrt{b_n}),
\]
where $\bar{D}_{n2}=\sum_{i=1}^n\{\mathbf{z}^{\top}_i[\E\mathbf{
S}_n(t_i)]^{-1}\mathbf{T}_n(t_i)\}^2$.
\end{lemma}
\begin{pf}
Note that
$D_{n2}-\bar{D}_{n2}=\sum_{i=1}^n\Gamma_1(i)\Gamma_2(i),$
where $\Gamma_1(i)=\mathbf{z}^{\top}_i(\mathbf{S}^{-1}_n(t_i)+[\E
\mathbf{
S}_n(t_i)]^{-1})\times \mathbf{T}_n(t_i)$ and $\Gamma_2(i)=\mathbf{
z}^{\top}_i(\mathbf{S}^{-1}_n(t_i)-[\E\mathbf
{S}_n(t_i)]^{-1})\mathbf{
T}_n(t_i)$.

Let $S \Gamma_1(i)=\sum_{j=1}^i\Gamma_1(i)$ for $1\le i \le n$ and
$S\Gamma_1(0)=0$. Then
\begin{eqnarray*}
D_{n2}-\bar{D}_{n2}=\sum_{i=1}^n
\bigl(S \Gamma_1(i)-S \Gamma_1(i)\bigr)
\Gamma_2(i)=\sum_{i=1}^{n-1}S
\Gamma_1(i) \bigl(\Gamma_2(i)-\Gamma_2(i+1)
\bigr)+S \Gamma_1(n)\Gamma_2(n).
\end{eqnarray*}
Note that
\begin{eqnarray*}
S \Gamma_1(i)&=&r_n^2\sum
_{k=1}^n\sum_{j=1}^i
\mathbf{z}^{\top}_j\bigl(\mathbf{ S}^{-1}_n(t_j)+
\bigl[\E\mathbf{S}_n(t_j)\bigr]^{-1}
\bigr)K_{b_n}(t_k-t_j)\pmatrix{\mathbf{x}_k \varepsilon_k
\vspace*{2pt}\cr
\mathbf{x}_k \varepsilon_k \bigl[(t_k-t_j)/b_n
\bigr]
}
\\
&=&r_n^2
\sum_{k=1}^n\Xi_1(i,k)
\varepsilon_k+r_n^2\sum
_{k=1}^n\Xi_2(i,k)\varepsilon_k,
\end{eqnarray*}
where
\begin{eqnarray*}
\Xi_1(i,k)&=&\sum_{j=1}^i
\mathbf{z}^{\top}_j\bigl(\mathbf{S}^{-1}_n(t_j)+
\bigl[\E\mathbf{ S}_n(t_j)\bigr]^{-1}
\bigr)K_{b_n}(t_k-t_j)\mathbf{z}^{\top}_k,
\\
\Xi_2(i,k)&=&\sum_{j=1}^i
\mathbf{z}^{\top}_j\bigl(\mathbf{S}^{-1}_n(t_j)+
\bigl[\E\mathbf{ S}_n(t_j)\bigr]^{-1}
\bigr)K_{b_n}(t_k-t_j) \bigl(
\mathbf{0}^{\top}_p,\mathbf{x}^{\top
}_k
\bigr)^{\top}.
\end{eqnarray*}
By Lemma \ref{lemsninv} and the H\"{o}lder's
inequality, ${\max_{i}}\|\Xi_1(i,k)\|\le Cnb_n$. Hence by similar
conditioning arguments as those in the proof
Lemma \ref{lembias1},
\[
r_n^2\max_i \Biggl\|\sum
_{k=1}^n\Xi_1(i,k)\varepsilon_k
\Biggr\|=\mathrm{O}(\sqrt{n}).
\]
Similarly, $r_n^2{\max_i}\|\sum_{k=1}^n\Xi_2(i,k)\varepsilon_k\|
=\mathrm{O}(\sqrt{n})$. Hence, ${\max_i}\|S \Gamma_1(i)\|=\mathrm{O}(\sqrt{n})$. By similar
arguments, we have
\[
\max_i\bigl\|\Gamma_2(i)-\Gamma_2(i+1)\bigr\|=\mathrm{O}
\bigl(r_n^4\bigr) \quad\mbox{and}\quad\bigl \|\Gamma_2(n)
\bigr\|=\mathrm{O}\bigl(r_n^2\bigr).
\]
Therefore
\begin{eqnarray*}
\E|D_{n2}-\bar{D}_{n2}|&\le&\sum
_{i=1}^{n-1}\bigl\|S \Gamma_1(i)\bigr\|\bigl\|
\Gamma_2(i)-\Gamma_2(i+1)\bigr\|+\bigl\|S \Gamma_1(n)\bigr\|
\bigl\|S\Gamma_2(n)\bigr\|
\\
&=&\mathrm{O}\bigl(1/\bigl(\sqrt{n}b^{2}_n\bigr)\bigr)=\mathrm{o}(1/
\sqrt{b_n}).
\end{eqnarray*}
The lemma follows.
\end{pf}

%Note that
%$$D_{n2}=\sum_{j=1}^n\sum_{i=1}^n\mathbf{z}^{\top}_i\mathbf{
%S}^{-1}_n(t_i)\mathbf{z}^*_{ij}K_{b_n}(t_j-t_i),$$ where $\mathbf{
%z}^*_{ij}=(\mathbf{x}^{\top}_j,\mathbf{x}^{\top}_j(t_j-t_i)/b_n)^{
%By Lemma \ref{lemsninv} and similar arguments as those in the proof
%Lemma \ref{lembias1}, we have
%Therefore $\E|D^2_{n2}-\bar{D}^2_{n2}|\le
%(\|D_{n2}\|+\|\bar{D}_{n2}\|)\|D_{n2}-\bar{D}_{n2}\|=\mathrm{O}(1/
%The lemma follows.

%le5 #&#
%
\begin{lemma}\label{lemasy3}
Under condition \textup{(A)} and the assumption that $nb_n^3\rightarrow\infty$,
we have
\[
\bar{D}_{n2}=\Theta_n+\mathrm{o}_\p(1/
\sqrt{b_n}),
\]
where $\Theta_n=\sum_{i=1}^n\mathbf{T}^{\top}_n(t_i)[\E\mathbf{
S}_n(t_i)]^{-1}\E[\mathbf{z}_i\mathbf{z}^{\top}_i][\E\mathbf{
S}_n(t_i)]^{-1}\mathbf{T}_n(t_i)$.
\end{lemma}
\begin{pf}
Note that
$\bar{D}_{n2}=\sum_{i=1}^n\mathbf{T}^{\top}_n(t_i)[\E\mathbf{
S}_n(t_i)]^{-1}\mathbf{z}_i\mathbf{z}^{\top}_i[\E\mathbf{
S}_n(t_i)]^{-1}\mathbf{T}_n(t_i).$ Therefore
\[
\bar{D}_{n2}-\Theta_n=\sum_{i=1}^n
\mathbf{T}^{\top}_n(t_i)\Theta_n(i),
\]
where
$\Theta_n(i)=[\E\mathbf{
S}_n(t_i)]^{-1}\{\mathbf{z}_i\mathbf{z}^{\top}_i-\E[\mathbf
{z}_i\mathbf{z}^{\top
}_i]\}[\E\mathbf{
S}_n(t_i)]^{-1}\mathbf{T}_n(t_i)$.
Note that
\begin{eqnarray*}
\sum_{j=1}^i\Theta_n(j)&=& r_n^2
\sum_{k=1}^n\sum
_{j=1}^i\bigl[\E\mathbf{ S}_n(t_j)
\bigr]^{-1}\bigl\{\mathbf{z}_j\mathbf{z}^{\top}_j-
\E\bigl[\mathbf{z}_j\mathbf{z}^{\top
}_j\bigr]
\bigr\}\bigl[\E\mathbf{ S}_n(t_j)\bigr]^{-1}\\
&&\hspace*{40pt}{}\times K_{b_n}(t_k-t_j)
\pmatrix{\mathbf{x}_k
\varepsilon_k
\vspace*{2pt}\cr
\mathbf{x}_k \varepsilon_k \bigl[(t_k-t_j)/b_n
\bigr]
}.
\end{eqnarray*}
By the short memory property of $\mathbf{x}_i$ in condition (A4) and
similar arguments as those in the proof of Lemma \ref{lembias1}, we
have
\begingroup
\abovedisplayskip=7pt
\belowdisplayskip=7pt
\begin{eqnarray*}
\max_i \Biggl\|\sum_{j=1}^i
\bigl[\E\mathbf{ S}_n(t_j)\bigr]^{-1}\bigl\{
\mathbf{z}_j\mathbf{z}^{\top}_j-\E\bigl[
\mathbf{z}_j\mathbf{z}^{\top
}_j\bigr]\bigr\}
\bigl[\E\mathbf{ S}_n(t_j)\bigr]^{-1}K_{b_n}(t_k-t_j)
\Biggr\|=\mathrm{O}(\sqrt{nb_n}).\vadjust{\goodbreak}
\end{eqnarray*}
Hence by similar conditioning arguments as those in the proof of Lemma
\ref{lembias1}, we have
\[
\max_{i}\Biggl\|\sum_{j=1}^i
\Theta_n(j)\Biggr\|=\mathrm{O}(\sqrt{n}r_n).
\]
Together with \eqref{eq4} and the summation by parts technique used in
Lemma \ref{lemasy1}, it follows that
$\E|\bar{D}_{n2}-\Theta_n|=\mathrm{O}(1/(\sqrt{n}b_n^{2}))=\mathrm{o}(1/\sqrt{b_n})$.
The lemma follows.
\end{pf}

%Under condition (A), we have
%where $\Theta_n=r_n^4\sum_{i=1}^n\sum_{j=1}^n\varepsilon_i\mathbf{
%z}^{\top}_i\E[\Xi_n(i,j)]\varepsilon_j\mathbf{z}_j$,
%$$\Xi_n(i,j)=\sum_{k=1}^n[\E\mathbf{S}_n(t_k)]^{-1}\mathbf{z}^*_{ik}(
%z}^*_{jk})^{\top}[\E\mathbf{
%S}_n(t_k)]^{-1}K_{b_n}(t_i-t_k)K_{b_n}(t_j-t_k),$$ and $\mathbf{
%z}^*_{ik}=((\mathbf{x}^{\top}_k,\mathbf{x}^{\top}_k(t_k-t_i)/b_n)^{
%Note that
%$$\bar{D}_{n2}=r_n^4\sum_{i=1}^n\sum_{j=1}^n\varepsilon_i\mathbf{
%z}^{\top}_i\Xi_n(i,j)\varepsilon_j\mathbf{z}_j.$$ Therefore
%z}^{\top}_i\Theta_n(i),
%where
%$\Theta_n(i)=\sum_{j=1}^n\{\Xi_n(i,j)-\E[\Xi_n(i,j)]\}\varepsilon_j
%z}_j$. By the short memory conditions of $(\varepsilon_i)$ and
%$(\mathbf{x}_i)$, we have
%and $\|\Theta_n(n)\|=\mathrm{O}(nb_n)$. Then by the summation by parts
%technique in Lemma \ref{lemasy2}, we have
%The lemma follows.

%Under conditions (A1)-(A7) and bandwidth conditions ($n^{-1/4}\ll
%b_n\le cn^{-2/9}$), we have uniformly on $[b_n,1-b_n]$,

%le6 #&#
%
\begin{lemma}\label{lemasy4}
Assume condition \textup{(A)}. Then on a possibly richer probability space,
there exist i.i.d standard $p$ dimensional Gaussian random vectors
$V_1,\ldots,V_n$, such that
%
%e30 #&#
%
\begin{equation}
\label{eq1} \bigl|\Theta_n-\Theta_n^*\bigr|+\bigl|\bar{D}_{n1}-
\bar{D}_{n1}^*\bigr|=\mathrm{O}_\p\bigl((\log n)^{3/2}/
\bigl(n^{1/4}b_n^{3/2}\bigr)\bigr),
\end{equation}
where
\begin{eqnarray*}
\Theta_n^*&=&\sum_{i=1}^n
\tilde{\mathbf{T}}^{\top}_n(t_i)\bigl[\E\mathbf{
S}_n(t_i)\bigr]^{-1}\E\bigl[
\mathbf{z}_i\mathbf{z}^{\top}_i\bigr] \bigl[\E
\mathbf{ S}_n(t_i)\bigr]^{-1}\tilde{
\mathbf{T}}_n(t_i),
\\[-2pt]
\bar{D}_{n1}^*&=&\sum
_{i=1}^n\tilde{V}^{\top}_i
\bigl[\E\mathbf{ S}_n(t_i)\bigr]^{-1}\tilde{
\mathbf{T}}_{n}(t_i).
\end{eqnarray*}
\end{lemma}
\begin{pf}
Recall the definitions of $\tilde{V}_i$, $\tilde{\mathbf{T}}_{n}(t)$ and
$\tilde{\mathbf{T}}_{n,l}(t)$ in Proposition \ref{propwildboots}. We will
only prove $\Theta_n-\Theta_n^*=\mathrm{O}_\p((\log
n)^{3/2}/(n^{1/4}b_n^{3/2}))$ since $\bar{D}_{n1}-\bar{D}_{n1}^*=\mathrm{O}_\p
((\log n)^{3/2}/(n^{1/4}b_n^{3/2}))$ follows by similar arguments. Note that
\begin{eqnarray*}
\Theta_n=\sum_{i=1}^n
\mathbf{T}^{\top}_n(t_i)\bigl[\E\mathbf{
S}_n(t_i)\bigr]^{-1}\E\bigl[
\mathbf{z}_i\mathbf{z}^{\top}_i\bigr] \bigl[\E
\mathbf{ S}_n(t_i)\bigr]^{-1}
\mathbf{T}_n(t_i):=\sum_{i=1}^n
\mathbf{T}^{\top}_n(t_i)\tilde{
\Theta}_n(i).
\end{eqnarray*}
By Corollaries 1 and 2 of Wu and Zhou \cite{WuZho}, on a possibly richer
probability space, there exist i.i.d $p$ dimensional standard Gaussian
random vectors $V_1,\ldots,V_n$, such that
%
%e31 #&#
%
\begin{equation}
\label{eq3} \max_{1\le i\le n}|\Delta_i|=\mathrm{O}_\p
\bigl(n^{1/4}(\log n)^{3/2}\bigr),
\end{equation}
where $\Delta_i=\sum_{j=1}^i(\varepsilon_j\mathbf{x}_j-\Lambda
^{1/2}(t_j)V_j)$. Write $\Theta^{(1)}_n=\sum_{i=1}^n\tilde{\mathbf{
T}}^{\top}_n(t_i)\tilde{\Theta}_n(i)$. Then
\begin{eqnarray*}
&&\bigl|\Theta_n-\Theta^{(1)}_n\bigr|\\[-2pt]
&&\quad=\Biggl |\sum
_{i=1}^n\bigl[\mathbf{T}^{\top
}_n(t_i)-
\tilde{\mathbf{T}}^{\top}_n(t_i)\bigr]\tilde{
\Theta}_n(i) \Biggr|
\\
&&\quad= \Biggl|\sum_{i=1}^n
\bigl[\bigl(\mathbf{T}^{\top}_{n,0}(t_i),
\mathbf{0}^{\top
}_p\bigr)-\bigl(\tilde{\mathbf{T}}^{\top}_{n,0}(t_i),
\mathbf{0}^{\top}_p\bigr)\bigr]\tilde{\Theta
}_n(i)+\bigl[\bigl(\mathbf{0}^{\top}_p,
\mathbf{T}^{\top}_{n,1}(t_i)\bigr)-\bigl(
\mathbf{0}^{\top
}_p,\tilde{\mathbf{T}}^{\top}_{n,1}(t_i)
\bigr)\bigr]\tilde{\Theta}_n(i)\Biggr |
\\
&&\hspace*{-3pt}\quad:= \Biggl|\sum
_{i=1}^n\bigl[W^{\top}_{n,0}(t_i)
\tilde{\Theta}_n(i)+W^{\top
}_{n,1}(t_i)
\tilde{\Theta}_n(i)\bigr] \Biggr|.
\end{eqnarray*}
\endgroup
Write $\tilde{\Delta}_i=(\Delta^{\top}_i, \mathbf{0}^{\top
}_p)^{\top}$ and
$\tilde{\Delta}_0=0$. Note that
\begin{eqnarray*}
\sum_{i=1}^nW^{\top}_{n,0}(t_i)
\tilde{\Theta}_n(i)&=&r_n^2\sum
_{i=1}^n\sum_{k=1}^n(
\tilde{\Delta}_k-\tilde{\Delta}_{k-1})K_{b_n}(t_k-t_i)
\tilde{\Theta}_n(i)
\\
&=&r_n^2\sum
_{k=1}^n(\tilde{\Delta}_k-\tilde{
\Delta}_{k-1})\sum_{i=1}^nK_{b_n}(t_k-t_i)
\tilde{\Theta}_n(i)\\
&:=&r_n^2\sum
_{k=1}^n(\tilde{\Delta}_k-\tilde{
\Delta}_{k-1})\Omega_n(k).
\end{eqnarray*}
By the summation by parts formula,
\begin{eqnarray*}
\Biggl|\sum_{k=1}^n(\tilde{
\Delta}_k-\tilde{\Delta}_{k-1})\Omega_n(k)\Biggr|&=&\Biggl|
\sum_{k=1}^{n-1}\tilde{\Delta}_k
\bigl(\Omega_n(k)-\Omega_n(k+1)\bigr)+\tilde{\Delta
}_n\Omega_n(n)\Biggr|
\\
&\le&\max_{1\le i\le n}|\tilde{
\Delta}_i|\Biggl(\sum_{k=1}^{n-1}\bigl|
\Omega_n(k)-\Omega_n(k+1)\bigr|+\bigl|\Omega_n(n)\bigr|
\Biggr).
\end{eqnarray*}
By the smoothness of $K(\cdot)$ and the similar arguments as those in
the proof of Lemma \ref{lembias1}, it follows that
\begin{eqnarray*}
\max_{1\le k\le n-1}\bigl\|\Omega_n(k)-\Omega_n(k+1)
\bigr\|=\mathrm{O}(r_n), \qquad\bigl\|\Omega_n(n)\bigr\|=\mathrm{O}(1/r_n).
\end{eqnarray*}
Therefore by \eqref{eq3}, we have
\[
\Biggl|\sum_{i=1}^nW^{\top}_{n,0}(t_i)
\tilde{\Theta}_n(i)\Biggr|=\mathrm{O}_\p\bigl\{n^{1/4}
\log^{3/2}n\bigl(nr^3_n+r_n\bigr)
\bigr\}=\mathrm{O}_\p\bigl((\log n)^{3/2}/\bigl(n^{1/4}b_n^{3/2}
\bigr)\bigr).
\]
Similarly,
\[
\Biggl|\sum_{i=1}^nW^{\top}_{n,1}(t_i)
\tilde{\Theta}_n(i)\Biggr|=\mathrm{O}_\p\bigl((\log n)^{3/2}/
\bigl(n^{1/4}b_n^{3/2}\bigr)\bigr).
\]
Hence, $|\Theta_n-\Theta^{(1)}_n|=\mathrm{O}_\p((\log
n)^{3/2}/(n^{1/4}b_n^{3/2}))$. Note that
\begin{eqnarray*}
\Biggl|\Theta^{(1)}_n-\sum_{i=1}^n
\tilde{\mathbf{T}}^{\top}_n(t_i)\bigl[\E\mathbf{
S}_n(t_i)\bigr]^{-1}\E\bigl[
\mathbf{z}_i\mathbf{z}^{\top}_i\bigr] \bigl[\E
\mathbf{ S}_n(t_i)\bigr]^{-1}\tilde{
\mathbf{T}}_n(t_i)\Biggr|=\sum_{i=1}^n
\hat{\Theta}_n(t_i)\bigl[\mathbf{T}_n(t_i)-
\tilde{\mathbf{T}}_n(t_i)\bigr],
\end{eqnarray*}
where $\hat{\Theta}_n(t_i)=\tilde{\mathbf{T}}^{\top}_n(t_i)[\E
\mathbf{
S}_n(t_i)]^{-1}\E[\mathbf{z}_i\mathbf{z}^{\top}_i][\E\mathbf{
S}_n(t_i)]^{-1}$.
Hence by similar arguments, it follows that
\[
\Biggl|\sum_{i=1}^n\hat{\Theta}_n(t_i)
\bigl[\tilde{\mathbf{T}}_n(t_i)-\mathbf{
T}_n(t_i)\bigr]\Biggr|=\mathrm{O}_\p\bigl((\log
n)^{3/2}/\bigl(n^{1/4}b_n^{3/2}\bigr)
\bigr).
\]
The lemma follows.
\end{pf}

%le7 #&#
%
\begin{lemma}\label{lemasy5}
Under condition \textup{(A)} and the assumption that $b_n\rightarrow0$,
$nb_n\rightarrow\infty$, we have
\begin{eqnarray*}
\sqrt{b_n}\biggl\{\Theta_n^*-2\bar{D}_{n1}^*-
\tilde{K}(0)\int_{0}^1\tr\bigl[H(t)H^{\top}(t)
\bigr] \,\mathrm{d}t/b_n\biggr\}\Rightarrow N\bigl(0,
\sigma^2\bigr).
\end{eqnarray*}
\end{lemma}
\begin{pf}
Note that both $\Theta^*_n$ and $D^*_{n1}$ are quadratic forms of
i.i.d. standard Gaussian random vectors. By Lemma \ref{lemsninv} and
similar arguments as those in the proof of Lemma \ref{lemasy3}, it can
be shown that $\Theta_n^*-\Theta^{**}_{n}=\mathrm{O}_\p(1)$ and $\bar
{D}_{n1}^*-\bar{D}_{n1}^{**}=\mathrm{O}_\p(1)$, where
\begin{eqnarray*}
\Theta^{**}_n&=&\sum_{i=1}^n
\tilde{\mathbf{T}}^{\top
}_{n,0}(t_i)M^{-1}(t_i)
\tilde{\mathbf{T}}_{n,0}(t_i),
\\
\bar{D}_{n1}^{**}&=&
\sum_{i=1}^n{V}^{\top}_i
\Lambda^{1/2}(t_i)M^{-1}(t_i)\tilde{
\mathbf{T}}_{n,0}(t_i).
\end{eqnarray*}
Note that
\begin{eqnarray*}
\Theta^{**}_n=r_n^4\sum
_{k=1}^n\sum_{r=1}^n
V_k^{\top}\Lambda^{1/2}(t_k)\Biggl[\sum
_{i=1}^nM^{-1}(t_i)K_{b_n}(t_k-t_i)K_{b_n}(t_r-t_i)
\Biggr]\Lambda^{1/2}(t_r)V_r
\end{eqnarray*}
and that $M^{-1}(t_i)K_{b_n}(t_k-t_i)K_{b_n}(t_r-t_i)=0$ if
$|t_k-t_r|\ge2b_n$ or $\min\{|t_i-t_r|,|t_i-t_k|\}\ge b_n$. Hence by
Lemma \ref{lemsninv} and similar arguments as those in the proof of
Lemma \ref{lemasy3}, it follows that
\[
\Theta_n^{**}-\Theta^{***}_{n}=\mathrm{O}(1)
\qquad\mbox{where } \Theta^{***}_{n}=r_n^2
\sum_{k=1}^n\sum
_{r=1}^n V_k^{\top}
\tilde{H}(t_k)K\ast K_{b_n}(t_k-t_r)
\tilde{H}^{\top}(t_r)V_r,
\]
where $\tilde{H}(\cdot)=\Lambda^{1/2}(\cdot)M^{-1/2}(\cdot)$. Similarly,
\[
D_{n1}^{**}-D_{n1}^{***}=\mathrm{O}(1)\qquad \mbox{where } D_{n1}^{***}=r_n^2\sum
_{k=1}^n\sum_{r=1}^n
V_k^{\top}\tilde{H}(t_k)K_{b_n}(t_k-t_r)
\tilde{H}^{\top}(t_r)V_r.
\]
Using the fact that $V_i$'s are i.i.d. standard Gaussian, elementary
calculations show that
\begin{eqnarray*}
\sqrt{b_n}\biggl\{\Theta_n^{***}-2
\bar{D}_{n1}^{***}-\tilde{K}(0)\int_{0}^1
\tr\bigl[H(t)\bigr] \,\mathrm{d}t/b_n\biggr\}\Rightarrow N\bigl(0,
\sigma^2\bigr).
\end{eqnarray*}
The lemma follows.
\end{pf}

%le8 #&#
%
\begin{lemma}\label{lemasy6}
Under conditions \textup{(A1)--(A7)}, we have
\[
\frac{\sum_{i=1}^n\varepsilon_i^2}{n}= \int_{0}^1
\vartheta^2(t) \,\mathrm{d}t+\mathrm{O}_\p(1/\sqrt{n}),
\]
where $\vartheta^2(t)=\E[ V(t, \FF_0)]^2$.
\end{lemma}
\begin{pf}
Note that $\E\varepsilon_i^2=\vartheta^2(t_i)$. Therefore
\[
\sum_{i=1}^n\bigl[\varepsilon_i^2-
\vartheta^2(t_i)\bigr]=\sum_{k=-\infty}^n
\sum_{i=1}^n\PP^*_k
\varepsilon_i^2,
\]
where $\PP^*_i(\cdot)=\E(\cdot|\RR_i)-\E(\cdot|\RR_{i-1})$. Since
$\PP^*_i$ and $\PP_{j}^*$ are orthogonal for $i\neq j$, we have
\begin{eqnarray*}
\Biggl\|\sum_{i=1}^n\bigl[\varepsilon_i^2-
\vartheta^2(t_i)\bigr]\Biggr\|^2=\sum
_{i=1}^n\sum_{j=1}^n
\sum_{k=-\infty}^n\E\bigl[\PP^*_k
\varepsilon_i^2\PP^*_k\varepsilon_j^2
\bigr]\le\sum_{i=1}^n\sum
_{j=1}^n\sum_{k=-\infty}^n
\bigl\|\PP^*_k\varepsilon_i^2\bigr\|\bigl\|
\PP^*_k\varepsilon_j^2\bigr\|.
\end{eqnarray*}
Let $(\chi^*_k)$ be an i.i.d. copy of $(\chi_k)$. By Theorem 1 in Wu
\cite{wu05}, $\|\PP^*_k\varepsilon_i^2\|\le
\|\varepsilon_i^2-\varepsilon_{i,k}^2\|$, where
$\varepsilon_{i,k}=(\RR_{k-1}, \chi^*_{k},\chi_{k+1},\ldots,
\chi_{i})$ if $k\le i$ and $\varepsilon_{i,k}=\varepsilon_{i}$
otherwise. By the Cauchy--Schwarz inequality, we have for $i\ge k$
\begin{eqnarray*}
&&\bigl\|\varepsilon_i^2-\varepsilon_{i,k}^2
\bigr\|\\
&&\quad\le\|\varepsilon_i+\varepsilon_{i,k}\|_4\|
\varepsilon_i-\varepsilon_{i,k}\|_4\le C
\bigl\|H(t_i,\GG_i)V(t_i,\FF_i)-H(t_i,
\GG_{i,k})V(t_i,\FF_{i,k})\bigr\|_4
\\
&&\quad\le C\bigl\{\bigl\|H(t_i,\GG_i)\bigr\|_4
\bigl\|(V(t_i,\FF_i)-V(t_i,\FF_{i,k})
\bigr\|_4+\bigl\|V(t_i,\FF_{i,k})\bigr\|_4
\bigl\|H(t_i,\GG_i)-H(t_i,\GG_{i,k})
\bigr\|_4\bigr\}
\\
&&\quad\le C(i-k+1)^{-2}.
\end{eqnarray*}
Therefore,
\begin{eqnarray*}
\Biggl\|\sum_{i=1}^n\bigl[\varepsilon_i^2-
\vartheta^2(t_i)\bigr]\Biggr\|^2\le C\sum
_{i=1}^n\sum_{j=1}^n
\sum_{k=-\infty}^{\min
(i,j)}(i-k+1)^{-2}(j-k+1)^{-2}
\le Cn.
\end{eqnarray*}
Hence,
$\|\sum_{i=1}^n[\varepsilon_i^2-\vartheta^2(t_i)]\|=\mathrm{O}(\sqrt{n})$.
Note that
$\sum_{i=1}^n\vartheta^2(t_i)=n\int_{0}^1\vartheta^2(t) \,\mathrm
{d}t+\mathrm{O}(1)$.
The lemma follows.
\end{pf}

%le9 #&#
%
\begin{lemma}\label{lemsninv}
Recall that $\mu_h=\int_{-1}^1x^hK(x) \,\mathrm{d}x$. Under
condition \textup{(A)}, we have
\begin{eqnarray*}
\sup_{0\le t\le1}\bigl\|\mathbf{S}_n^{-1}(t)-\bigl[\E\mathbf{
S}_n(t)\bigr]^{-1}\bigr\|_{8}=\mathrm{O} \biggl(
\frac{1}{\sqrt{nb_n}} \biggr).
\end{eqnarray*}
Additionally, $\sup_{0\le t\le1}|[\E\mathbf{S}_n(t)]^{-1}|=\mathrm{O}(1)$. For
$h=0,2$, we have
\begin{eqnarray*}
\sup_{b_n \le t \le1-b_n} \bigl|\bigl[\E\mathbf{S}_{n, h}(t)\bigr]^{-1}
- \bigl[\mu_h M(t)\bigr]^{-1}\bigr| = \mathrm{O}\bigl(b_n^2
\bigr).
\end{eqnarray*}

\end{lemma}
\begin{pf}
The proof follows by the similar arguments as those of Lemma 6 in
Zhou and Wu \cite{ZhoWu10}. Details are omitted.
\end{pf}

\begin{pf*}{Proof of Theorem \ref{thm1}}
Theorem \ref{thm1} follows from Lemmas \ref{lembias1}--\ref
{lemasy6} above and the Slutsky's theorem.
\end{pf*}

\begin{pf*}{Proof of Theorem \ref{thm2}}
Recall that $\varepsilon_i=([\varepsilon^{(1)}_{i}]^{\top
},[\varepsilon^{(2)}_{i}]^\top)^{\top}$ and
$V_i=([V^{(1)}_{i}]^{\top
},[V^{(2)}_{i}]^\top)^{\top}$, where $\varepsilon^{(1)}_{i}$ and
$V^{(1)}_{i}$ are $p_1$ dimensional. Note that, under $H_{01}$, we have
a local linear regression of $y^*_i$ on $\mathbf{x}_i^{(2)}$. Recall again
the definitions of $\mathbf{S}^{(2)}_{n}$, $\mathbf{S}^{(2)}_{nl}$,
$\mathbf{
z}_{i}^{(2)}$, $\tilde{V}^{(2)}_i$, $\mathbf{T}^{(2)}_{n}$, $\mathbf{
T}^{(2)}_{nl}$, $\tilde{\mathbf{T}}^{(2)}_{n}$ and $\tilde{\mathbf{
T}}^{(2)}_{nl}$ in Section~\ref{secwildboots}.

Following very similar arguments as those in Lemmas \ref{lembias1} to
\ref{lemasy6}, it can be shown that
%
%e32 #&#
%
\begin{equation}
\label{eq5} \RSS_1-\RSS_0-\mathbf{B}^{(2)}_n=
\Theta^{(2)*}_n-2\bar{D}_{n1}^{(2)*}+\mathrm{o}_{\p
}(1/
\sqrt{b_n}),
\end{equation}
where $\mathbf{B}^{(2)}_n=\frac{nb_n^4\mu_2^{2}}{4}\int_{0}^1\{
[\betav^{(2)}(t)]''\}^\top M_{22}(t)[\betav^{(2)}(t)]'' \,\mathrm
{d}t+\mathrm{o}_\p(nb_n^4)$,
\begin{eqnarray*}
\Theta_n^{(2)*}&=&\sum_{i=1}^n
\bigl[\tilde{\mathbf{T}}^{(2)}_n(t_i)
\bigr]^{\top}\bigl[\E\mathbf{ S}^{(2)}_n(t_i)
\bigr]^{-1}\E\bigl[\mathbf{z}^{(2)}_i
\mathbf{z}^{(2)\top}_i\bigr] \bigl[\E\mathbf{
S}^{(2)}_n(t_i)\bigr]^{-1}\tilde{
\mathbf{T}}^{(2)}_n(t_i),
\\
\bar{D}_{n1}^{(2)*}&=&\sum_{i=1}^n
\tilde{V}^{(2)\top}_i\bigl[\E\mathbf{ S}^{(2)}_n(t_i)
\bigr]^{-1}\tilde{\mathbf{T}}^{(2)}_{n}(t_i).
\end{eqnarray*}
Note that $\Theta_n^{(2)*}$ and $\bar{D}_{n1}^{(2)*}$ are quadratic
forms of i.i.d. Gaussian vectors $V_1,\ldots,V_n$. Theorem \ref{thm2}
follows easily from \eqref{eq1} and \eqref{eq5}.
\end{pf*}

\begin{pf*}{Proof of Proposition \ref{propequasemi}}
Define $\mathbf{Y}^*=(y^*_1,\ldots,y^*_n)^{\top}$ and $\tilde
{\mathbf{
Y}}^*=(\tilde{y}^*_1,\ldots,\tilde{y}^*_n)^{\top}$. Let $\hat
{\varepsilon}_i$ and $\tilde{\varepsilon}_i$ be the $i$th residual of
the local linear regression of $y^*_i$ and $\tilde{y}^*_i$ on $\mathbf{
x}_{i}^{(2)}$, respectively. From \eqref{eqsol}, we can write
$\hat{\varepsilon}_i=y^*_i-R_i\mathbf{Y}^* \mbox{ and } \tilde
{\varepsilon
}_i=\tilde{y}^*_i-R_i\tilde{\mathbf{Y}}^*$,
where $R_i$ is a $1\times n$ vector which can be written in a closed
form \eqref{eqsol}. Note also that $R_i$ is functionally independent
of the errors $\varepsilon_i$. Hence,
\begin{eqnarray*}
\widetilde{\RSS}_1-\RSS_1=\sum
_{i=1}^n\bigl(\tilde{\varepsilon}^2_i-
\hat{\varepsilon}^2_i\bigr)=\sum
_{i=1}^n(\tilde{\varepsilon}_i-\hat{
\varepsilon}_i)^2+2\sum_{i=1}^n
\hat{\varepsilon}_i(\tilde{\varepsilon}_i-\hat{
\varepsilon}_i):=I+2\mathit{II}.
\end{eqnarray*}
Let $\delta_i=-(\mathbf{x}^{(1)}_i)^{\top}(\betav
^{(1)}_{0}(t_i,\hat{\theta
})-\betav^{(1)}_{0}(t_i,\theta_0))$ and $\Delta_n=(\delta_1,\ldots
,\delta_n)$. Hence,
\begin{eqnarray*}
\E(I)=\sum_{i=1}^n\|\tilde{
\varepsilon}_i-\hat{\varepsilon}_i\|^2=\sum
_{i=1}^n\|\delta_i-R_i
\Delta_n\|^2.
\end{eqnarray*}
From condition (B), it is easy to see that, for sufficiently large $n$,
\begin{eqnarray*}
\Bigl\|\max_{1\le i\le n}\bigl|\betav^{(1)}_{0}(t_i,
\hat{\theta})-\betav^{(1)}_{0}(t_i,
\theta_0)\bigr|\Bigr\|_4=\mathrm{O}(1/\sqrt{n}).
\end{eqnarray*}
Therefore, it is easy to derive from condition (A) that
%
%e33 #&#
%
\begin{eqnarray}
\label{eq7} \max_{1\le i \le n}\|\delta_i\|=\mathrm{O}(1/\sqrt{n})\quad\mbox{and}\quad \max_{1\le i
\le n}\|R_i\Delta_n\|=\mathrm{O}(1/\sqrt{n}).
\end{eqnarray}
Hence, $I=\mathrm{O}_\p(1)$. We now deal with $\mathit{II}$. Note that, by \eqref{eq6},
\begin{eqnarray*}
\hat{\varepsilon}_i=\varepsilon_i-\bigl(
\mathbf{z}^{(2)}_i\bigr)^{\top}\bigl(\hat{\etav
}^{(2)}(t_i)-\etav^{(2)}(t_i)\bigr)=
\varepsilon_i-\bigl(\mathbf{z}^{(2)}_i
\bigr)^{\top} \bigl(\mathbf{S}_n^{(2)}(t_i)
\bigr)^{-1}\bigl[\mathbf{B}^{(2)}_n(t_i)+
\mathbf{T}^{(2)}_n(t_i)\bigr].
\end{eqnarray*}
Hence,
\begin{eqnarray*}
\mathit{II}&=&\sum_{i=1}^n\bigl(
\mathbf{z}^{(2)}_i\bigr)^{\top} \bigl(
\mathbf{S}_n^{(2)}(t_i)\bigr)^{-1}
\mathbf{B}^{(2)}_n(t_i)[\tilde{\varepsilon
}_i-\hat{\varepsilon}_i]+\sum
_{i=1}^n\varepsilon_i[\tilde{
\varepsilon}_i-\hat{\varepsilon}_i]
\\
&&{}+\sum
_{i=1}^n\bigl(\mathbf{z}^{(2)}_i
\bigr)^{\top} \bigl(\mathbf{S}_n^{(2)}(t_i)
\bigr)^{-1}\mathbf{T}^{(2)}_n(t_i)[
\tilde{\varepsilon}_i-\hat{\varepsilon}_i]
\\
&:=&\mathit{II}^{*}+\mathit{II}^{**}+\mathit{II}^{***}.
\end{eqnarray*}
By H\"{o}lder inequality, condition (A) and \eqref{eq7}, the bias term
\[
\E\bigl|\mathit{II}^{*}\bigr|\le\sum_{i=1}^n\bigl\|
\mathbf{z}^{(2)}_i\bigr\|_6\bigl\|\mathbf{
S}_n^{(2)}(t_i)^{-1}\bigr\|_6
\bigl\|\mathbf{B}^{(2)}_n(t_i)\bigr\|_6\bigl[\|
\delta_i\|+\| R_i\Delta_n\|\bigr]=\mathrm{O}\bigl(
\sqrt{n}b_n^2\bigr).
\]
Write $J_i=-(\mathbf{x}^{(1)}_i)^{\top}\frac{\partial\betav
^{(1)}_{0}(t_i,\theta_0)}{\partial\theta}$ and let $\mathbf
{J}=(J^{\top
}_1,\ldots,J^{\top}_n)^{\top}$. By second order Taylor expansion of
$\betav^{(1)}_{0}(t_i,\hat{\theta})$ at $\theta_0$ and condition (B),
it is easy to see that
%
%e34 #&#
%
\begin{equation}
\tilde{\varepsilon}_i-\hat{\varepsilon}_i=(J_i-R_i
\mathbf{J}) (\hat{\theta}-\theta_0)+r_i,
\end{equation}
with the reminder term $r_i$ satisfying $\max_{1\le i \le n}\|r_i\|
=\mathrm{O}(1/n)$. Therefore,
\begin{eqnarray*}
\E\bigl|\mathit{II}^{**}\bigr|\le\Biggl\|\sum_{i=1}^n
\varepsilon_i(J_i-R_i\mathbf{J})\Biggr\|\bigl\|(\hat{
\theta}-\theta_0)\bigr\|+\max_{1\le i \le n}\|r_i\|\sum
_{i=1}^n\|\varepsilon_i\|
\end{eqnarray*}
By the similar conditioning arguments as those in the proof of Lemma
\ref{lembias1}, it is easy to show that $\|\sum_{i=1}^n\varepsilon
_i(J_i-R_i\mathbf{J})\|=\mathrm{O}(\sqrt{n})$. Hence $\E|\mathit{II}^{**}|=\mathrm{O}(1)$. By
similar arguments and elementary but tedious calculations, it follows
that $\E|\mathit{II}^{***}|=\mathrm{O}(1)$. Therefore, the proposition follows.
\end{pf*}

\begin{pf*}{Proof of Proposition \ref{proplocalpower}}
Let $\overline{\RSS}_0=\sum_{i=1}^n\varepsilon_i^2$. Then $\RSS
_a-\RSS_0=\RSS_a-\overline{\RSS}_0- (\RSS_0-\overline{\RSS}_0)$.
Under the
local alternative $\betav(\cdot)=\betav_0(\cdot)+n^{-4/9}\mathbf
{f}_n(\cdot
)$, we have
\begin{eqnarray*}
\RSS_0-\overline{\RSS}_0=n^{-4/9}\sum
_{i=1}^n\mathbf{f}_n^{\top}(t_i)
\mathbf{ x}_i\varepsilon_i+n^{-8/9}\sum
_{i=1}^n\bigl[\mathbf{f}_n^{\top}(t_i)
\mathbf{x}_i\bigr]^2.
\end{eqnarray*}
By the similar arguments as those in the proof of Lemma \ref
{lembias1}, it is easy to show that
\begin{eqnarray*}
\sqrt{b_n}n^{-8/9}\sum_{i=1}^n
\bigl[\mathbf{f}_n^{\top}(t_i)\mathbf{
x}_i\bigr]^2&=&c^{1/2}\int_0^1
\mathbf{f}_n^{\top}(t)M(t)\mathbf{f}_n(t) \,
\mathrm{d}t+\mathrm{o}_\p(1),
\\
\sum_{i=1}^n
\mathbf{f}_n^{\top}(t_i)\mathbf{x}_i
\varepsilon_i&=&\mathrm{O}_\p\bigl(n^{1/2}\bigr).
\end{eqnarray*}
On the other hand, by Lemmas \ref{lembias1}--\ref{lemasy6} and the
fact that $\betav(\cdot)=\betav_0(\cdot)+n^{-4/9}\mathbf
{f}_n(\cdot)$, it
is easy to show that
\begin{eqnarray*}
\hspace*{-4pt}&&\sqrt{b_n}\biggl\{\RSS_a-\overline{\RSS}_0-
\frac{\tilde{K}(0)}{b_n}\int_{0}^1\tr\bigl[H(t)\bigr]
\,\mathrm{d}t\biggr\}-\frac{c^{9/2}\mu_2^{2}}{4}\int_0^1
\bigl[\betav''(t)\bigr]^\top M(t)
\betav''(t) \,\mathrm{d}t
-\frac{c^{9/2}\mu_2^{2}}{4}F_2
\\
\hspace*{-4pt}&&\quad\Rightarrow N\bigl(0,\sigma^2\bigr)
\end{eqnarray*}
and $\RSS_0/n=\mathcal{ V}+\mathrm{o}_\p(1)$. Therefore, the proposition follows.
\end{pf*}

\begin{pf*}{Proof of Theorem \ref{thmnewtest}}
A careful check of Lemmas \ref{lembias1} and \ref{lembias2} shows
that the asymptotic bias of $\lambda_n^*$
%
%e35 #&#
%
\begin{eqnarray}
\label{eqnewbias} \mathbf{B}_n^*=\int_{c_{\min}}^{c_{\max}}
\frac{n(zn^{-\gamma})^4\mu_2^{2}}{4} \,\mathrm{d}z \int_0^1\bigl[
\betav''(t)\bigr]^\top M(t)
\betav''(t) \,\mathrm{d}t+\mathrm{o}_\p
\bigl(n^{1-4\gamma}\bigr).
\end{eqnarray}
Another careful check of Lemmas \ref{lemasy1} to \ref{lemasy6} and
using Lemma \ref{lemsninv} show that
%
%e36 #&#
%
\begin{eqnarray}
\label{eqnewvar} \lambda_n^*-\mathbf{B}_n^*&=&\int
_{c_{\min}}^{c_{\max}}\sum_{k=1}^n
\sum_{r=1}^nV_k^{\top}
\tilde{H}(t_k)\bigl[2K_{zn^{-\gamma}}(t_k-t_r)
\nonumber
\\[-8pt]
\\[-8pt]
\nonumber
&&\hspace*{56pt}{}-K\ast K_{zn^{-\gamma}}(t_k-t_r)\bigr]
\tilde{H}^{\top
}(t_r)V_r/\bigl(nzn^{-\gamma}
\bigr) \,\mathrm{d}z+\mathrm{o}_\p\bigl(n^{-\gamma/2}\bigr).\qquad
\end{eqnarray}
Since $V_i$'s are i.i.d. standard Gaussian, a central limit theorem for
$\lambda_n^*-\mathbf{B}_n^*$ can be easily derived. Now Theorem \ref
{thmnewtest} follows from \eqref{eqnewbias} and \eqref{eqnewvar}.
Details are omitted.
\end{pf*}

\begin{pf*}{Proof of Proposition \ref{proppower}}
By the Cauchy--Schwarz inequality,
\begin{eqnarray*}
\int_\R Q(c_{max},y)^2 \,
\mathrm{d}y&=&\int_\R\biggl[\int_{c_{\min}}^{c_{\max
}}
\bigl(\bigl[2K(y/z)-K\ast K(y/z)\bigr]/\sqrt{z}\bigr)\times
1/\sqrt{z} \,
\mathrm{d}z \biggr]^2 \,\mathrm{d}y
\\
&\le& \int_\R
\biggl[\int_{c_{\min}}^{c_{\max}}\bigl[2K(y/z)-K\ast K(y/z)
\bigr]^2/z \,\mathrm{d}z\int_{c_{\min}}^{c_{\max}}1/z
\,\mathrm{d}z\biggr] \,\mathrm{d}y
\\
&=&\bigl(\log(c_{\max})-
\log(c_{\min})\bigr)\int_{c_{\min}}^{c_{\max}}\int
_\R\bigl[2K(y/z)-K\ast K(y/z)\bigr]^2/z \,
\mathrm{d}y \,\mathrm{d}z
\\
&=&\bigl(\log(c_{\max})-
\log(c_{\min})\bigr) (c_{\max}-c_{\min})\int
_{\R}\tilde{K}^2(t) \,\mathrm{d}t.
\end{eqnarray*}
Consider any fixed $c\in(0,\infty)$. Plugging the above inequality into
\eqref{powernewtest} and letting $c_{\max}\downarrow c$ and $c_{\min
}\uparrow c$, it follows that $\sup_{0<c_{\min}<c_{\max}<\infty
}\beta^*_{\alpha}(c_{\min},c_{\max})\ge\beta_{\alpha}(c)$.
Hence, the proposition follows.
\end{pf*}

\section*{Acknowledgements}
I am grateful to the two anonymous
referees for their many helpful
comments which greatly improved the quality of the original version of
the paper. The research was supported in part by NSERC of Canada.

% imsref loaded by akundreckaite, 2012-12-13 09:52:54
%

\printhistory

\end{document}